\documentclass[journal,twoside,web]{ieeecolor}
\usepackage{generic}
\usepackage{cite}
\usepackage{amsmath,amssymb,amsfonts}
\usepackage{algorithmic}
\usepackage{graphicx}
\usepackage{algorithm,algorithmic}
\usepackage{hyperref}
\usepackage{textcomp}
\usepackage{amsmath}
\usepackage{diagram}
\usepackage[normalem]{ulem}
\usepackage{bm}
\usepackage{mathrsfs}
\usepackage[utf8]{inputenc}

\allowdisplaybreaks[1]
\makeatletter
\newenvironment{breakablealgorithm}
{% \begin{breakablealgorithm}
  \begin{center}
   \refstepcounter{algorithm}% New algorithm
   \hrule height.8pt depth0pt \kern2pt% \@fs@pre for \@fs@ruled
   \renewcommand{\caption}[2][\relax]{% Make a new \caption
    {\raggedright\textbf{\ALG@name~\thealgorithm} ##2\par}%
    \ifx\relax##1\relax % #1 is \relax
    \addcontentsline{loa}{algorithm}{\protect\numberline{\thealgorithm}##2}%
    \else % #1 is not \relax
    \addcontentsline{loa}{algorithm}{\protect\numberline{\thealgorithm}##1}%
    \fi
    \kern2pt\hrule\kern2pt
   }
  }{% \end{breakablealgorithm}
  \kern2pt\hrule\relax% \@fs@post for \@fs@ruled
 \end{center}
}
\makeatother

\hypersetup{hidelinks=true}

\newtheorem{note}{Remark}[section]
\newtheorem{assumption}{Assumption}[section]
\newtheorem{lemma}{Lemma}[section]
\newtheorem{proposition}{Proposition}[section]
\newtheorem{theorem}{Theorem}[section]

\def\BibTeX{{\rm B\kern-.05em{\sc i\kern-.025em b}\kern-.08em
    T\kern-.1667em\lower.7ex\hbox{E}\kern-.125emX}}
\begin{document}
\title{A Unified Recursive Identification Algorithm with Quantized Observations Based on Weighted Least-Squares Type Criteria}
\author{Xingrui Liu, \IEEEmembership{Student Member, IEEE}, Ying Wang, \IEEEmembership{Member, IEEE}, and  Yanlong Zhao, \IEEEmembership{Senior Member, IEEE}
\thanks{ 
This work is supported by the National Natural Science Foundation of China under Grants 62025306, 62588101, 62303452, and T2293773, CAS Project for Young Scientists in Basic Research under Grant YSBR-008, China Postdoctoral Program for Innovation Talents under BX20230403. (Corresponding author: Yanlong Zhao.)}
\thanks{Xingrui Liu, Ying Wang, and Yanlong Zhao are with the State Key Laboratory of Mathematical Sciences, Academy of Mathematics and Systems Science, Chinese Academy of Sciences, Beijing 100190, P. R. China.
Xingrui Liu and Yanlong Zhao are also with the School of Mathematics Sciences, University of Chinese Academy of Sciences, Beijing 100149, P. R. China.
Ying Wang is also with the Division of Decision and Control Systems, KTH Royal Institute of Technology, Stockholm 11428, Sweden (e-mail: liuxingrui@amss.ac.cn; wangying96@amss.ac.cn; ylzhao@amss.ac.cn).
}
\thanks{This is the accepted version of the paper to appear in IEEE Transactions on Automatic Control, doi: 10.1109/TAC.2025.3620623.}
}

\pubid{%
  \makebox[\textwidth][c]{% 在一栏宽度内居中
    \parbox{\textwidth}{\centering\footnotesize
© 2025 IEEE.  Personal use of this material is permitted.   Permission from IEEE must be obtained for all other uses, in any current or future media, including reprinting/republishing this material for advertising or promotional purposes, creating new collective works, for resale or redistribution to servers or lists, or reuse of any copyrighted component of this work in other works.
    }%
  }%
}

\maketitle

\begin{abstract}
This paper investigates system identification problems with Gaussian inputs and quantized observations under fixed thresholds.
By reinterpreting the nonlinear effects induced by quantization as the product of the unknown parameter and an unknown nonlinear coefficient, this work establishes a novel weighted least-squares criterion that enables linear estimation of unknown parameters under quantized observations.
Subsequently, a two-step recursive identification algorithm is constructed by estimating two unknown terms, which is capable of handling both Gaussian noisy and noise-free linear systems.
Convergence analysis of this identification algorithm is conducted, demonstrating convergence in both almost sure and $L^{p}$ senses under mild conditions, with respective rates of $O(\sqrt{ \log \log k/k})$ and $O(1/k^{p/2})$, where $k$ denotes the time step.
In particular, this algorithm offers an asymptotically efficient estimation of the variance of Gaussian variables using quantized observations.
Furthermore, extensions to output-error systems are discussed, enhancing the applicability and relevance of the proposed methods.
Two numerical examples are provided to validate these theoretical advancements.
\end{abstract}


\section{Introduction} \label{sec:introduction}

\IEEEPARstart{D}{ue} to the cost-effectiveness and power limitations inherent in sensor technologies, the challenge of quantized identification has become increasingly prominent across various systems, including industrial systems \cite{anwar2004anti}, \cite{sun2004aftertreatment}, \cite{cai2021dissipative}, networked systems \cite{jiang2024linear}, and even biological systems \cite{ghysen2003origin}.
Quantized identification denotes estimating the unknown parameters of systems based on quantized measurements, where only the discrete set to which the system output belongs is discernible, without exact knowledge of the output value.
Over the past two decades, numerous studies have emerged addressing quantized identification challenges, proposing various methods including empirical measure method  \cite{zhang2003system}, \cite{yin2007asymptotically}, \cite{guo2015asymptotically}, expectation maximization method \cite{godoy2011identification}, \cite{marelli2013identification}, \cite{zhao2016iterative}, stochastic approximation method \cite{guo2013recursive},  \cite{song2018recursive}, \cite{huang2024identification}, \cite{wang2023identification}, stochastic gradient method \cite{zhang2022identification}, \cite{wang2023asymptotically}, \cite{wang2024threshold} and sign-error method  \cite{csaji2012recursive}, \cite{zhao2017recursive},  \cite{you2015recursive}, \cite{wang2022unified}.
Indeed, the prevailing algorithms are primarily grounded on the following three fundamental techniques:

The first technique is to utilize the noise distributions.  
For instance, \cite{zhang2003system} derived a probability expression for binary-valued observations using a noise distribution.
They developed an identification algorithm under periodic input conditions employing the empirical measure method based on the inverse function of the noise distribution.
Similarly, under general bounded persistently exciting inputs, \cite{guo2013recursive} utilized noise distributions to establish a probabilistic relationship between quantized observations and parameters, leading to the development of a stochastic approximation type recursive estimator for finite impulse response (FIR) systems. 
Subsequent studies have expanded upon these foundations, addressing issues such as more generalized system models \cite{wang2023identification}, weaker excitations \cite{guo2015asymptotically, zhang2022identification}, and enhanced convergence properties \cite{yin2007asymptotically, wang2023asymptotically}.

\pubidadjcol

The second technique revolves around adaptive thresholds, which are specially designed to compare the size relationship between actual and predicted outputs. 
For instance, under independent and identically distributed (iid) inputs,  \cite{csaji2012recursive} and \cite{zhao2017recursive}  investigated quantized identification problems of ARX systems and nonlinear FIR systems, respectively, utilizing the sign-error type algorithm with expanding truncations proposed by \cite{chen2003asymptotic}.
Reference \cite{you2015recursive} introduced a stochastic approximation type identification algorithm for linear systems with adaptive binary-valued observations.
Furthermore, these adaptive thresholds make it possible to estimate the unknown parameter in both noisy and noise-free cases.
Under general bounded persistently exciting inputs, \cite{wang2022unified} presented a unified sign-error type identification algorithm for FIR systems.

The third technique entails the design of tailored inputs to enhance parameter information within finite quantized data.
For example, \cite{casini2007time} and \cite{casini2011input} addressed system identification problems in a worst-case setting with binary-valued observations.
Employing two distinct cost functions: maximum parametric identification error and relative uncertainty reduction concerning the minimum achievable error, they devised suboptimal input signals for identifying FIR models and provided an upper bound for time complexity.

However, it is essential to note that the first technique relies on noise distributions, rendering it inadequate for handling identification problems in noise-free cases. 
The second technique encounters challenges in some practical systems where the time-invariant quantizers with fixed thresholds cannot be designed \cite{wang2023identification}.
Moreover, designing thresholds in the second technique and designing inputs in the third technique may entail significant costs.
These techniques are subject to fundamental limitations.
In the cases where thresholds are fixed and inputs cannot be designed, almost all existing works have yet to consider developing recursive algorithms that simultaneously adapt to both noisy and noise-free scenarios.

Motivated by these considerations, this study investigates quantized identification problems under fixed thresholds, encompassing scenarios both with and without noise. 
Specifically, this paper utilizes optimization methods.

In classical system identification theory, optimization methods such as the weighted least-squares (WLS) algorithm have been proven effective in resolving system identification problems in both noisy and noise-free cases \cite{ljung1987theory}. 
However, when only quantized observations are available, the absence of accurate output information presents a considerable challenge in computing prediction errors.
Consequently, the classical WLS algorithm cannot be directly applied.

Relevant studies have proposed criteria aimed at addressing quantized identification challenges by substituting prediction errors in the classical WLS with those derived from quantized observations \cite{735317}, \cite{colinet2009weighted}, \cite{jafari2012convergence}.
These criteria introduced a predictor of quantized observations by directly applying the quantizer to the predicted accurate outputs.
The primary challenge in developing a recursive identification algorithm based on these criteria stems from the nonlinear effects induced by quantization, which complicates the derivation of differential expressions for the predicted quantized observations concerning the parameter estimations.
To address this issue, \cite{735317} utilized the approximation of the derivative of the quantizer.
Reference \cite{colinet2009weighted} incorporated predicted outputs into the weight coefficients and subsequently proposed an optimization criterion.
Based on this criterion, \cite{jafari2012convergence} treated the predicted quantized observations as an independent variable in parameter estimations and developed a recursive identification algorithm accordingly.
However, based on these criteria, the identification algorithms mentioned above can only operate offline or theoretically be proved convergent in noise-free cases. 

Recently, \cite{song2018recursive} proposed a stochastic approximation quantized identification algorithm by establishing a correlation function between quantized observations, Gaussian inputs, and system parameters, which exhibits almost sure convergence. 
Building on this work, \cite{mestrah2023identification} introduced a least-squares type quantized identification algorithm proven to converge in the mean square sense.
These algorithms apply to both Gaussian noise and noise-free conditions. 
However, due to the presence of product terms involving unknown parameters and their nonlinear functions in the correlation functions, these algorithms require prior knowledge of the parameters.
For instance, the algorithm in \cite{song2018recursive} is limited to the semiparametric response model, while \cite{mestrah2023identification} can only identify the parameter direction.
Therefore, based on these works, this paper further investigates criterion formulation and algorithm development.
The primary contributions of this paper are as follows:

\romannumeral1) This paper proposes a novel WLS-type criterion for quantized identification by replacing the prediction errors of accurate observations with those of quantized observations.
Previous methods \cite{735317}, \cite{colinet2009weighted}, \cite{jafari2012convergence} predict quantized observations by directly applying the quantizer to predicted accurate outputs, which yields quantized predictions that are nonlinear and non-differentiable concerning the system parameter, complicating the design of recursive algorithms.
In contrast, we skillfully reinterpret the nonlinear effects induced by quantization as the product of the unknown parameter and an unknown nonlinear coefficient, and thereby introduce a novel quantized observation predictor that retains linearity in the parameter estimation, facilitating a recursive identification algorithm by estimating the two unknown terms.

\romannumeral2) This paper develops a two-step unified recursive identification algorithm suitable for linear systems in both Gaussian noise and noise-free scenarios.
Firstly, we design a maximum likelihood type algorithm to estimate the unknown nonlinear coefficient, which is represented by the system output variance.
Based on the constructed nonlinear coefficient estimates, we derive the system parameter estimates by minimizing the WLS-type criterion.
This identification algorithm applies to both cases since the variance estimates and the WLS-type criterion are designed to handle both Gaussian noise and noise-free cases.
Additionally, it leverages richer statistical information from Gaussian inputs rather than deterministic ones, eliminating the need for designable thresholds as in \cite{wang2022unified}. 

\romannumeral3) This paper demonstrates that the algorithm converges in both almost sure and $L^{p}$ senses with rates of \( O(\sqrt{\log \log k / k}) \) and \( O(1 / k^{p/2}) \), respectively, matching the fast convergence rates attained by the classical WLS algorithm under accurate observations \cite{ljung1987theory, radenkovic2000almost}. 
Additionally, the algorithm provides an asymptotically efficient estimation of the variance of Gaussian variables using quantized observations in the first step.
Furthermore, it removes the need for prior knowledge of the semiparametric response model as in \cite{song2018recursive, huang2024identification} or the parameter norms as in \cite{mestrah2023identification}, since the second step of the algorithm estimates the parameter direction while jointly estimating their norm with the first step.

\romannumeral4) 
This paper further tackles the quantized identification problem for output-error (OE) systems, where the system output exhibits infinite impulse response dynamics with respect to the system inputs.
By applying Durbin’s method from \cite{durbin1960fitting}, we transform the OE system identification problem into a dynamic regression task and solve it using the proposed two-step algorithm. 
To address the core non-iid challenge in dynamic regression convergence analysis, we construct auxiliary processes that are strictly stationary and geometrically strongly mixing, along with controllable error terms to characterize the non-iid sequences.
Based on this construction, we prove that, the proposed algorithm achieves an almost sure convergence rate of $O(\sqrt{\log \log k / k})$, matching the classical WLS rate under accurate observations \cite{radenkovic2000almost}, and improving upon existing $O(1/k^{v})$ rate for any $v \in (0, 1/2)$ in prior works \cite{song2018recursive}, \cite{huang2024identification}, \cite{wang2023identification}.

The remainder of this paper is organized as follows. 
Section \ref{sec_a} formulates the problem with linear system structures and discusses the system identifiability.
Section \ref{sec_b} focuses on the algorithm construction and establishes the main results.
Section \ref{sec_f} further extends the methods and results to OE systems.
Section \ref{sec_c} presents the proofs of the main results.
Section \ref{sec_d} gives simulation examples to verify the conclusion.
Section \ref{sec_g} is the summary and prospect of this paper.

$\mathbf{Notations:}$
In this paper, $\mathbb{R}$, $\mathbb{Z}$, $\mathbb{R}^{n}$  and $\mathbb{R}^{m \times n}$ are the sets of real number, integers, $n$-dimensional real vectors, and matrices with rows $m$ and columns $n$, respectively. 
For a constant $x$, $\vert x \vert$ denotes its absolute value.
For a pair of integers $N \geq k \geq 0$,  ${N}\choose{k}$ denotes its binomial coefficients.
For a vector $a=[a_{1},a_{2},\ldots ,a_{n}]^{T} \in \mathbb{R}^{n}$, 
$\Vert a \Vert$ denotes its Euclidean norm, i.e, $\Vert a \Vert = (\sum_{i=1}^{n} a_{i}^{2})^{1/2}$.
$\mathbf{1}_{n} = [1,1,\ldots,1]^{T} \in \mathbb{R}^{n}$.
$\mathbf{0}_{n} = [0,0,\ldots,0]^{T} \in \mathbb{R}^{n}$.
For a matrix $A$, $A^T$ denotes its transpose; $A_{i,j}$ denotes its element in the $i$-th row and the $j$-th column; $\text{rank}(A)$ denotes its rank; 
$A^{-1}$ denotes its inverse matrix.
$I_{n}$ is an $n$-dimension identity matrix.
$\mathbb{P}$ denotes the probability operator.
$\mathbb{E}$ denotes the expectation operator.
$I_{\{\cdot\}}$ denotes the indicator function, whose value is $1$ if its argument (a formula) is true and $0$ otherwise.
$F(\cdot)$ and $f(\cdot)$ denote the cumulative distribution and probability density functions of the standard Gaussian random variable.
$F^{-1}(\cdot)$ is the inverse function of $F(\cdot)$. 
$\mathcal{N}(\mu,\delta^2)$ denotes the (multivariate) Gaussian distribution with mean $\mu$ and standard deviation $\delta$.

\section{Problem Formulation} \label{sec_a}

Consider a linear system described by
\begin{align}\label{model_a}
y_{k} = \phi_{k}^{T}\theta + d_{k},\quad k \geq 1,
\end{align}
where $k$ is the time index; $\theta \in \mathbb{R}^{n}$ is a vector of unknown parameters; $\phi_{k} \in \mathbb{R}^{n}$ is the system input; $d_{k} \in \mathbb{R}$ is the potential system noise. Especially if $d_{k} = 0, k \geq 1$, the linear system (\ref{model_a}) is precisely a noise-free system. The system output $y_{k}$ cannot be exactly measured and can only be measured by a quantized observation:
\begin{align}\label{model_b}
s_{k} = Q(y_{k}) =
\left\{
\begin{array}{lcl}
0, & \text{if} & y_{k}  \leq  C_{1}, \\
1, & \text{if} & C_{1} < y_{k} \leq C_{2}, \\
\vdots & & \vdots \\
m, & \text{if} & y_{k} > C_{m},
\end{array}   
\right.
\end{align}
where $Q(\cdot)$ is the quantizer; $-\infty = C_{0} < C_{1}  <  C_{2} <\ldots< C_{m} < C_{m+1} = \infty$ are the known thresholds; $m$ is the number of the quantizer thresholds. The quantized observation can also be represented as $s_{k} = \sum_{i=0}^{m}i I_{\{ C_{i} < y_{k} \leq C_{i+1} \} }$.

The goal is to develop a unified algorithm to estimate the unknown parameter vector $\theta$ based on the system input $\{\phi_{k}\}_{k=1}^{\infty}$ and the quantized observation $\{ s_{k} \}_{k=1}^{\infty}$ in both noisy and noise-free cases.

\begin{note}
This paper starts with a static linear model to illustrate the algorithm construction. The extensions to high-dimensional OE systems are given in Section \ref{sec_f}.
\end{note}

\subsection{Assumptions}

To proceed with our analysis, we introduce some assumptions regarding the inputs and noise.

\begin{assumption}\label{ass_a}
(Stochastic persistent excitation)
$\{ \phi_{k} \}_{k=1}^{\infty}$ is a sequence of independent and identically distributed (iid) Gaussian random variables with zero mean and an unknown positive definite covariance matrix $H \in \mathbb{R}^{n \times n}$, i.e., $\phi_k \sim \mathcal{N}(\mathbf{0}_{n}, H)$, where $H \triangleq \mathbb{E}[\phi_{k}\phi_{k}^{T}] > 0, k \geq 1$.
\end{assumption}

\begin{note} 
The Gaussian input assumption has been widely applied across various fields, including system identification \cite{song2018recursive, colinet2009weighted, 6746194}, signal processing \cite{jacobsson2019linear}, and stochastic control \cite{grancharova2008explicit}.  
For instance, in the field of system identification, it is reasonable in a practical scenario, such as the fluctuations in pollutant concentrations in rivers \cite{9280402}.
Furthermore, it has been recognized as an effective excitation signal \cite{6650183}.  
In particular, properties such as the independence of uncorrelated jointly Gaussian variables and the light-tailed nature of their distributions are critical for the development of the unified recursive identification algorithm proposed in this work.
\end{note}

\begin{note} 
Under Assumption \ref{ass_a}, the system regressors in the linear system (\ref{model_a}) coincide with the system inputs and are therefore iid.
The more general case involving non-iid regressors in the high-dimensional OE systems will be addressed in Section \ref{sec_f}.
\end{note}

\begin{assumption}\label{ass_b}
(Potential Gaussian noise)
$\{d_{k}\}_{k=1}^{\infty}$ is a sequence of iid Gaussian random variables with zero mean and an unknown variance $\delta_{d}^{2}$, i.e., $d_{k} \sim \mathcal{N}(0, \delta_d^{2})$,  where $ \delta_{d}^{2} \triangleq  \mathbb{E}[ d_{k}^{2} ] \geq 0$ for all $k \geq 1$, and is independent of $\{ \phi_{k} \}_{k=1}^{\infty}$.
\end{assumption}

\begin{note}
This paper investigates the quantized identification problem under both Gaussian noise and noise-free conditions, allowing the noise variance $\delta_{d}^{2}$ to be either positive or zero.
Notably, this work does not require prior knowledge of the noise variance $\delta_{d}^{2}$, 
in contrast to the popular estimation methods, including empirical measure method \cite{zhang2003system}, approximate message passing method \cite{rangan2011generalized}, and maximum likelihood (ML) method \cite{10416234}, which require noise distributions as priors.
\end{note}

\subsection{System identifiability}

Due to the lack of information about accurate output $y_{k}$, one may discuss the system identifiability at first.

\begin{proposition}\label{pri_the_c}
If Assumptions \ref{ass_a} and \ref{ass_b} hold, then the linear system (\ref{model_a})-(\ref{model_b}) is unidentifiable with $m=1$ and $C_{1}=0$.
\end{proposition}

The proof of Proposition \ref{pri_the_c} is supplied in Section \ref{j8k}.

\begin{note}
In fact, in the cases where $m=1$ and $C_{1}=0$, quantized observation $s_{k}$ reflects the symbol information of system output $y_{k}$. 
In the absence of noise, equal scale amplification or reduction of parameters will not affect the symbol information of the output since it merely scales the system output.  
Similarly, in the presence of Gaussian noise with zero mean, such operations do not probabilistically affect the output symbol information due to the symmetry of the Gaussian noise in distribution.
Indeed, in this case, the linear system (\ref{model_a})-(\ref{model_b}) can be determined up to a scaling factor \cite{bottegal2017new}. 
\end{note}

Therefore, this paper considers all cases except for one mentioned in Proposition \ref{pri_the_c}.

\section{Algorithm design and its properties}\label{sec_b}

This section will construct a two-step recursive identification algorithm based on quantized observations in both noisy and noise-free cases and establish its convergence properties.

\subsection{The design concept of the identification algorithm}

We would like to introduce the design concept of the identification algorithm in this subsection first.

In scenarios where accurate output measurements are available, the WLS algorithm proves to be an effective identification method, applicable in both noisy and noise-free environments \cite{ljung1987theory}. 
The core idea of the WLS algorithm is to minimize the criterion \( \sum_{l=1}^{k} \beta_l (y_l - \phi_l^T \hat{\theta}_k)^2 \) to obtain the parameter estimate \( \hat{\theta}_k \), where \( \{\beta_l\}_{l=1}^{k} \) denotes weight coefficient.

However, when only quantized observations are accessible, obtaining the prediction error $y_{l} - \phi^{T}_{l}\hat{\theta}_{k}$ becomes unfeasible.
Furthermore, the nonlinearity introduced by the quantizer $Q(\cdot)$ poses challenges in providing the differential expression for the quantized observation $s_{l} = Q(y_{l})= Q(\phi^{T}_{l}\theta + d_{l})$ with respect to the system parameter $\theta$. 
Consequently, the majority of existing recursive quantized identification algorithms have not been formulated using optimization methods. 

To address this challenge, inspired by \cite{song2018recursive}, we establish a probabilistic relationship among quantized observation, Gaussian input, and system parameter as a preliminary step.

\begin{proposition} \label{pri_the_b}
If Assumptions \ref{ass_a} and \ref{ass_b} hold, then,
\begin{align}\label{pri_2}
\mathbb{E}\left[s_{k}\phi_{k}\right] = \rho\left(\delta_{y}\right) H \theta,
\end{align}
where $\rho(\delta_{y}) \triangleq \sum_{i=1}^{m}  \exp(-C_{i}^{2}/(2\delta_{y}^{2})) /(\sqrt{2 \pi} \delta_{y})$; 
$\delta_{y}$ denotes the standard deviation of $y_{k}$, i.e., $\delta_{y}^{2}  \triangleq \mathbb{E}[y_{1}^{2}] = \theta^{T}H\theta + \delta_{d}^{2}$. 
\end{proposition}

The proof of Proposition \ref{pri_the_b} is supplied in Section \ref{j7k}.

\begin{note}
Proposition \ref{pri_the_b} shows that the nonlinear effects induced by quantization $Q(\cdot)$ can be reinterpreted as the product of an unknown parameter $\theta$ and an unknown nonlinear coefficient $\rho(\delta_y)\phi_k$.
Specifically, since $\mathbb{E}[\phi_{k}\phi_{k}^{T}] = H$, the quatized observation $s_k$ can be reinterpreted as $\rho(\delta_y)\phi_k^T\theta$ plus a uncorrelated residual term. 
\end{note}

Based on this  Proposition \ref{pri_the_b}, we designate $\rho(\delta_{y})\phi_{l}^{T}\hat{\theta}_{k}$ as the predictor of the quantized observation $s_{l}$, and proposed the WLS-type criterion for quantized idnetification:
\begin{align}\label{ujijo2}
J_{k} = \sum\limits_{l=1}^{k} \beta_{l}\left( s_{l} - \rho\left(\delta_{y}\right)\phi_{l}^{T}\hat{\theta}_{k} \right)^{2},
\end{align}
where the prediction errors $y_{l} - \phi^{T}_{l}\hat{\theta}_{k}$ in the classical WLS criterion is substituted by the prediction errors of quantized observations $s_{l} - \rho(\delta_{y})\phi_{l}^{T}\hat{\theta}_{k}$.

\begin{note}
The weight coefficient $\beta_{k}$ enables the adjustment of the individual weight of prediction errors to emphasize different observations. For instance, assigning smaller weight coefficients to outliers can help reduce their impact
\cite{ljung1987theory}.
\end{note}

It's worth noting that $\delta_{y}^{2} = \theta^{T}H\theta + \delta_{d}^{2}$.
Consequently, $\rho(\delta_{y})$ in the WLS-type criterion (\ref{ujijo2}) depends on the unknown parameter $\theta$, making direct computation unattainable.

\begin{figure}[h]
\centering
\includegraphics[scale=0.24]{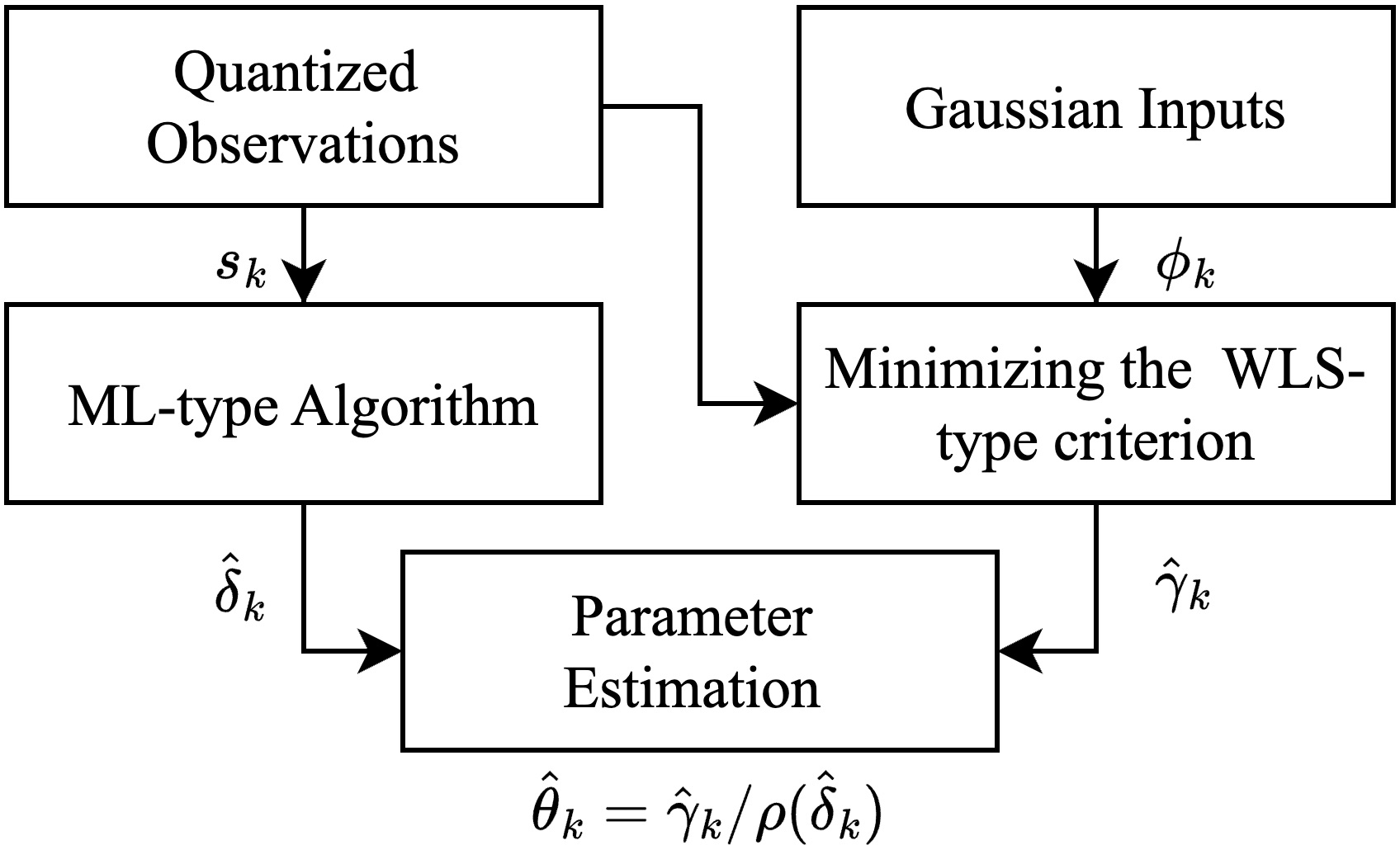}
\caption{Identification algorithm block diagram}
\label{fig_1}
\end{figure}

To address this challenge, we introduce $\gamma \triangleq \rho(\delta_y)\theta$ and propose a two-step algorithm: first, estimating $\delta_{y}$, and then estimating $\gamma$ to subsequently obtain $\theta$.
In Subsection \ref{sedu}, we introduce the ML-type algorithm for estimating $\delta_y$. 
Then, in Subsection \ref{oip}, we reformulate the WLS-type criterion (\ref{ujijo2}) as $J_{k} = \sum_{l=1}^{k} \beta_{l} ( s_{l} - \phi_{l}^{T}\hat{\gamma}_{k})^{2}$ and derive the optimal estimate of $\gamma$ under this criterion. 
Finally, the estimate of the unknown parameter $\theta$ is obtained as $\hat{\theta}_{k} = \hat{\gamma}_{k} / \rho(\hat{\delta}_{k})$, where $\hat{\delta}_k$ represents the estimate of $\delta_y$ and $\hat{\gamma}_k$ represents the estimate of $\gamma$ at time $k$. 
A simplified depiction of the design concept of the identification algorithm is illustrated in Fig. \ref{fig_1}. 

\begin{note}
An intuitive explanation of the two-step algorithm is as follows: the estimate of $\rho(\delta_y)$ from the first step corresponds to the norm of the unknown parameter, while the direction is determined by $\hat{\gamma}_k$ from the second step. 
Specifically, the second step uses information from the input $\phi_k$ to establish the direction of the estimate. 
In contrast, the first step treats the problem as an estimation of the variance of Gaussian variables based on quantized observations, which is independent of input information.
\end{note}

\subsection{First step: estimating the system output's variance} \label{sedu}

In this subsection, we will introduce the design concept and establish the convergence properties of the ML-type algorithm for estimating $\delta_y$. We treat this problem as estimating the variance of Gaussian variables based on quantized observations, independent of system identification. Given that the output $y_{k}$ is iid, the log-likelihood function can be expressed as:
\begin{align*}
l_{k} = \log \mathbb{P}\left(s_{1},s_{2},\ldots ,s_{k}|\delta_{y}\right) = \sum\limits_{l=1}^{k} \log \mathbb{P}\left(s_{l}|\delta_{y}\right).
\end{align*}
Define $S_{k}^{i}  \triangleq \sum_{l=1}^{k} I_{\{ s_{l} = i \}}/k, i = 0,1,2,\ldots,m$. Then by $\mathbb{P}(s_{l} = i|\delta_{y}) = \mathbb{P}(C_{i} < y_{l} \leq C_{i+1}|\delta_{y}) = F(C_{i+1}/\delta_{y})-F(C_{i}/\delta_{y})$, we have 
\begin{align}\label{ttr}
l_{k} =  k\sum\limits_{i=0}^{m} S_{k}^{i} \log\left(F\left(\frac{C_{i+1}}{\delta_{y}}\right)-F\left(\frac{C_{i}}{\delta_{y}}\right)\right).
\end{align}
It follows that
\begin{align}\label{tttu}
\frac{ \partial l_{k} }{ \partial \delta_{y} } =  k\sum\limits_{i=0}^{m} S_{k}^{i} \frac{ f(C_{i}/\delta_{y})C_{i}/\delta_{y}^{2}- f(C_{i+1}/\delta_{y})C_{i+1}/\delta_{y}^{2} }{ F(C_{i+1}/\delta_{y}) - F(C_{i}/\delta_{y}) }.
\end{align}

Specially, in the case where $m=1$, since $C_{1} \neq 0$, the solution of $\partial l_{k}/\partial \delta_{y}= 0$ yields the maximum likelihood estimation (MLE) of $\delta_y$ in this cases as:\footnote{Note that $F^{-1}(\cdot)$ is not invertible at $S_{k}^{0}=0$ or $S_{k}^{0}=1$ and $C_{1}/F^{-1}(\cdot)$  is not invertible at $S_{k}^{0}=1/2$. Based on the idea of EM method without truncation \cite{yin2007asymptotically}, one may modify these points by  $S_{k}^{0}=c^{*}$ when $S_{k}^{0}=0$, $S_{k}^{0}=1/2$ or $S_{k}^{0}=1$, where $c^{*} \in (0,1)$ and $c^{*} \neq 1/2$, which does not affect the convergence analysis and properties. Therefore, this modification will not be explicitly stated in the subsequent proofs and development. Besides, the non-invertible points with general quantized observations can be modified in the same way.}
\begin{align}\label{equ_mm1}
\hat{\delta}_{k} = \frac{C_{1}}{F^{-1}\left(S_{k}^{0}\right)}. 
\end{align}

In the case where $m>1$, as (\ref{tttu}) represents a nonlinear function, obtaining the explicit solution for the MLE of $\delta_y$ through differentiation becomes challenging. 
Inspired by (\ref{equ_mm1}), we address parameter estimation under binary-valued observations by considering different thresholds $C_{j}, j=1,2, \ldots, m$ separately, and then combine them to estimate $\delta_y$.

Specifically, for nonzero threshold $C_{j}$, the $j$-th estimate of $\delta_{y}$ is defined as $\hat{\delta}_{k}^{j}  = C_{j}/F^{-1}(\sum_{i=0}^{j-1} S_{k}^{i})$.
For zero threshold $C_{j}$, we utilize another threshold to estimate $\delta_y$: if $j > 1$, since $C_1 \neq 0$ and $\sum_{i=0}^{j-1} S_{k}^{i}$ is precisely the value of the $k$-sample empirical distribution of output $y_{k}$ at $C_{j}/\delta_y$, we define $\hat{\delta}_{k}^{j}=C_{1}/F^{-1}(F(0) - \sum_{i=1}^{j-1} S_{k}^{i})$. Similarly, if $C_{1} = 0$, then we define $\hat{\delta}_{k}^{1}=C_{m}/F^{-1}(F(0) + \sum_{i=j}^{m-1} S_{k}^{i})$. 
Based on these estimates, the estimate of $\delta_y$ at time $k$ could be given by
\begin{align}\label{equ_op}
\hat{\delta}_{k} = \sum_{j=1}^{m} \mu_j \hat{\delta}_{k}^{j} = \mu^T \hat{\Delta}_{k},
\end{align}
where $\mu = [\mu_1,\mu_2,\ldots ,\mu_m]^T\in \mathbb{R}^{m}; \mu_{j} \in \mathbb{R}, j = 1,2,\ldots,m$ satisfies $\sum_{j=1}^{m} \mu_j = 1$; $\hat{\Delta}_{k}=[\hat{\delta}_{k}^{1}, \hat{\delta}_{k}^{2}, \ldots , \hat{\delta}_{k}^{m}]^{T}\in \mathbb{R}^{m}$.

In cases where $m = 1$, (\ref{equ_op}) equals to the MLE (\ref{equ_mm1}). 
The variance of the estimation errors of the MLE asymptotically reaches the minimum variance of unbiased estimates, i.e., Cram{\' e}r–Rao (CR) lower bound \cite{van2000asymptotic}.
 Given that the variance of estimation errors can serve as a criterion to assess the effectiveness of an unbiased estimator, we aim to design $\mu$ in (\ref{equ_op}) to approach the performance of the MLE from the perspective of the CR lower bound.

Define $\Delta_y  \triangleq [\delta_y,\delta_y,\ldots ,\delta_y]^{T} \in \mathbb{R}^{m}$, $\tilde{\Delta}_{k}  \triangleq \hat{\Delta}_{k} - \Delta_{y}$ and $V_{k}  \triangleq \mathbb{E}[\tilde{\Delta}_{k}\tilde{\Delta}_{k}^{T}]$. Then, the variance of the estimation errors of (\ref{equ_op}) is given by
\begin{align}\label{hiujik12}
\mathbb{E}\left[\left(\hat{\delta}_{k}- \delta_{y}\right)^2\right] = \mu^T \mathbb{E}\left[\tilde{\Delta}_{k}\tilde{\Delta}_{k}^{T}\right]  \mu =  \mu^T V_{k} \mu.
\end{align}

Using (\ref{hiujik12}) as the criterion, the problem of designing $\mu$ can be transformed into the following optimization problem:
\begin{align}\label{w78tyi}
\begin{array}{ll}
\text{minmize} & \mu^T V_{k} \mu \\
\text{subject to} & \mu^{T}\mathbf{1}_{m} = 1.
\end{array}
\end{align}

Based on the Gauss–Markov estimation \cite{luenberger1997optimization}, the solution of the problem (\ref{w78tyi}) could be given by
\[
\mu^* = \frac{V_{k}^{-1}\mathbf{1}_{m}}{\mathbf{1}_{m}^{T}V_{k}^{-1}\mathbf{1}_{m}}.
\]

However, since $\delta_y$ is unknown, $V_{k}$ is not explicitly derived. 
Inspired by \cite{yin2007asymptotically}, we present a key proposition for estimating $V_{k}$.

\begin{proposition} \label{pri_the_oi} 
Under Assumptions \ref{ass_a} and \ref{ass_b},
\begin{align}\label{pri_3}
\lim\limits_{k \to \infty} kV_{k} =   \left(U+G^{T}\right)^{-1}\left(W-ww^T\right)\left(U+G\right)^{-1},
\end{align}
where $F_{i} = F(C_{i}/\delta_{y})$, $f_i = f(C_{i}/\delta_{y})C_{i}/\delta_y^2$,  $i = 1,2,$ $\ldots,$ $m$; $U = \text{diag}(f_{1},f_{2},\ldots ,f_m)\in \mathbb{R}^{m \times m}$; $w = [F_1,F_{2},\ldots ,F_m]^T\in \mathbb{R}^{m}$; $G \in \mathbb{R}^{m \times m}$ satisfies 
$$
G_{u,v} = 
\left\{
\begin{array}{ll}
f_1,  & \text{if} \, C_{j}=0, j >1, u=1,v=j, \\
-f_1, & \text{if}  \, C_{j}=0, j >1, u=j, v=j, \\
f_m, & \text{if}  \, C_{1}=0, u=m, v=1, \\
-f_m & \text{if}  \, C_{1}=0, u=1, v=1, \\
0, & \text{otherwise};
\end{array} 
\right. 
$$
$$
W =
 \begin{bmatrix}
F_{1} &F_{1} &\cdots  &F_{1}\\
F_{1}&F_{2}&\cdots  &F_{2} \\
\vdots & \vdots &\ddots &  \vdots \\
F_{1}& F_{2}&\cdots  &F_{m}
\end{bmatrix}.
$$
\end{proposition}

The proof of Proposition \ref{pri_the_oi} is supplied in Section \ref{iu2}.

Then, the ML-type algorithm for estimating $\delta_y$ in a recursive form is constructed as Algorithm \ref{alg_b}.
\begin{breakablealgorithm}
\begin{algorithmic}\label{alg_b}
\caption{The ML-type algorithm for estimating $\delta_y$}
Beginning with initial values  $\hat{\delta}_0 \in \mathbb{R}$ and $S_{0}^{i} = 0$, for $i = 0,1,2,\ldots,m-1$, $j = 1,2,\ldots, m$, the algorithm for estimating $\delta_y$ is recursively defined at any $k \geq 1$ as follows:
\begin{align}
\label{www23456}  \hat{\delta}_{k} = & \Pi\left(\left[\hat{\delta}_{k}^{1}, \hat{\delta}_{k}^{2}, \ldots ,\hat{\delta}_{k}^{m}\right]\hat{\mu}_{k}\right),\\
\label{ytrt2} \hat{\delta}_{k}^{j} = & \frac{C_{j}I_{\{C_{j} \neq 0\}} }{F^{-1}(\sum_{i=0}^{j-1} S_{k}^{i})}+ \frac{C_{1}I_{\{C_{j} = 0, j > 1 \}}}{F^{-1}(F(0) - \sum_{i=1}^{j-1} S_{k}^{i})} \nonumber\\
& + \frac{C_{m}I_{\{C_{j} = 0, j = 1 \}}}{F^{-1}(F(0)  + \sum_{i=j}^{m-1} S_{k}^{i})}, \\
S_{k}^{i} = & \frac{k-1}{k} S_{k-1}^{i}+\frac{1}{k}I_{\{ s_{k} = i \}}, \\
\label{www23451}  \hat{\mu}_{k} = & \frac{ (\hat{U}_{k}+\hat{G}_{k})(\hat{W}_{k}-\hat{w}_{k}\hat{w}_{k}^T)^{-1}(\hat{U}_{k}+\hat{G}_{k}^{T})\mathbf{1}_{m}}{\mathbf{1}_{m}^{T}(\hat{U}_{k}+\hat{G}_{k})(\hat{W}_{k}-\hat{w}_{k}\hat{w}_{k}^T)^{-1}(\hat{U}_{k}+\hat{G}_{k}^{T})\mathbf{1}_{m}},
\end{align}
where $\hat{\delta}_{k}$ is the estimate of $\delta_{y}$ at time $k$;
$\Pi(\cdot)$ denotes a projection operator as $\Pi(x) =  \text{argmin} \{ \Vert x - \zeta \Vert, \zeta \in [c,1/c], \forall x \in R \}$ with a sufficiently small positive constant $c$ for ensuring that $\hat{\mu}_{k}$ has a solution;
$\hat{U}_{k} = \text{diag}(\hat{f}_{k}^{1},\hat{f}_{k}^{2},\ldots ,\hat{f}_{k}^{m});$  $\hat{w}_{k} = [\hat{F}_{k}^1,\hat{F}_{k}^{2},\ldots ,\hat{F}_{k}^m]^T$; $\hat{G}_{k} \in \mathbb{R}^{m \times m}$ satisfies
$$(\hat{G}_{k})_{u,v} =
\left\{
\begin{array}{ll}
\hat{f}_{k}^{1},  &  \text{if} \, C_{j}=0, j >1, u=1,v=j, \\
-\hat{f}_{k}^{1}, &  \text{if} \, C_{j}=0, j >1, u=j, v=j, \\
\hat{f}_{k}^{m}, &  \text{if} \, C_{1}=0, u=m, v=1 \\
-\hat{f}_{k}^{m} &  \text{if} \, C_{1}=0, u=v=1, \\
0, & \text{otherwise};
\end{array} 
\right. 
$$
$$ \hat{W}_{k} =
\begin{bmatrix}
\hat{F}_{k}^{1} &\hat{F}_{k}^{1} &\cdots  &\hat{F}_{k}^{1}\\
\hat{F}_{k}^{1}&\hat{F}_{k}^{2}&\cdots  &\hat{F}_{k}^{2} \\
\vdots & \vdots & \ddots &  \vdots \\
\hat{F}_{k}^{1}& \hat{F}_{k}^{2}&\cdots  &\hat{F}_{k}^{m}
\end{bmatrix};
$$
$\hat{F}_{k}^{i} = F(C_{i}/\hat{\delta}_{k-1})$; $\hat{f}_{k}^{i} = f(C_{i}/\hat{\delta}_{k-1})C_{i}/(\hat{\delta}_{k-1})^{2}$.
\end{algorithmic}
\end{breakablealgorithm}

The ML-type algorithm has the following convergence properties and asymptotic efficiency.

\begin{theorem}\label{mai_the_b}
Under Assumptions \ref{ass_a} and \ref{ass_b}, $\hat{\delta}_{k}$ given by the ML-type algorithm has the following properties:

\romannumeral1) Almost sure convergence: $\hat{\delta}_{k}$ converges to $\delta_{y}$ in the almost sure sense with a convergence rate of $O(\sqrt{\log \log k/k})$, i.e.,
\begin{align}
\label{y7tgyghuh}& \left\vert \tilde{\delta}_{k} \right\vert = O\left(\sqrt{\frac{\log\log k}{k}}\right), \quad \text{a.s},
\end{align}
where $\tilde{\delta}_{k} \triangleq  \hat{\delta}_{k} - \delta_{y}$ is the estimation error of $\delta_{y}$.

\romannumeral2) $L^{p}$ convergence: $\hat{\delta}_{k}$ converges to $\delta_{y}$ in the $L^{p}$ sense with a convergence rate of $O(1/k^{p/2})$, i.e., 
\begin{align}
\label{the_1_2}& \mathbb{E}\left[\left\vert\tilde{\delta}_{k}\right\vert^{p} \right] = O\left(\frac{1}{k^{p/2}}\right),
\end{align}
where $p$ is an arbitrary positive integer.

\romannumeral3) Asymptotic efficiency: $\hat{\delta}_{k}$ is an asymptotically efficient estimate of $\delta_y$ based on quantized observations, i.e.,
\begin{align}
\label{theorem_3}\lim\limits_{k \to \infty} \mathbb{E}\left[\left\vert \sqrt{k}\tilde{\delta}_{k}\right\vert^{2} \right]-k\sigma_{\text{CR}}(k) = 0.
\end{align}
where $\sigma_{\text{CR}}(k) = (\sum_{i=0}^{m}\tilde{f}_{i}^2/\tilde{F}_{i})^{-1}/k$ is the CR lower bound for estimating the variance of Gaussian variables based on quantized observations at time $k$;  $\tilde{f}_{l} = f_{l+1} - f_{l}$, $\tilde{F}_{l} = F_{l+1} - F_{l}$, $l = 0,1,2,\ldots,m$.
\end{theorem}

The proof of Theorem \ref{mai_the_b} is supplied in Section \ref{sectiondeltay}.

\begin{note}
Analyzing convergence in \( L^{p} \) spaces for varying values of \( p \) offers a detailed framework for estimation error.
For instance, \( L^1 \) convergence guarantees asymptotic unbiasedness of the remote estimator,  \( L^2 \) convergence aligns with classical mean square convergence, and  \( L^4 \) convergence enhances robustness against outliers.  
\end{note}

\subsection{Second step: estimating the system parameter} \label{oip}
This subsection will introduce the design concept and establish the convergence properties of the WLS-type algorithm for estimating $\theta$.

Note that $\gamma = \rho(\delta_y)\theta$ and $\hat{\gamma}_k$ denotes the estimate of $\gamma$ at time $k$. The WLS-type criterion (\ref{ujijo2}) can be expressed as:
\begin{align*}
J_{k} = \sum\limits_{l=1}^{k}\beta_{l}\left( s_{l} - \phi_{l}^{T}\hat{\gamma}_{k} \right)^{2}.
\end{align*}
Following the derivation method of the classical WLS algorithm, we set $\partial J_{k}/\partial \hat{\gamma}_{k}^{*} = 0$, yielding $\hat{\gamma}_{k}^{*} = (\sum_{l=1}^{k} \beta_{l}\phi_{l}\phi_{l}^{T})^{-1}( \sum_{l=1}^{k} \beta_{l}s_{l}\phi_{l})$.
Then, similarly to the recursive form of the WLS algorithm, the WLS-type algorithm for estimating $\theta$ is constructed as Algorithm \ref{alg_a}.

\begin{breakablealgorithm}
\begin{algorithmic}\label{alg_a}
\caption{The WLS-type algorithm for estimating $\theta$}
Beginning with initial values $\hat{\gamma}_{0} \in \mathbb{R}$ and a positive definitive matrix $P_{0} \in \mathbb{R}^{n \times n}$, based on $\hat{\delta}_{k}$ by the ML-type algorithm, the algorithm for estimating $\theta$ is recursively defined at any $k \geq 1$ as follows:
\begin{align}
&  \label{theta22009} \hat{\theta}_{k} = \frac{\hat{\gamma}_{k}}{\rho(\hat{\delta}_{k})}, \\
& \hat{\gamma}_{k} = \hat{\gamma}_{k-1} + \alpha_{k}P_{k-1}\phi_{k}\left(s_{k} - \phi_{k}^{T}\hat{\gamma}_{k-1}\right),\\
& \alpha_{k} = \frac{1}{ \beta_{k}^{-1} + \phi_{k}^{T}P_{k-1}\phi_{k}}, \\
& P_{k} = P_{k-1} - \alpha_{k}P_{k-1}\phi_{k}\phi_{k}^{T}P_{k-1}, \\
& \label{rho989} \rho(\hat{\delta}_{k})  = \frac{1}{\sqrt{2 \pi} \hat{\delta}_{k}}\sum_{i=1}^{m}  \exp\left(-\frac{C_{i}^{2}}{2( \hat{\delta}_{k})^{2}}\right),
\end{align}
where $ \hat{\theta}_{k}$ is the estimate of $\theta$ at time $k$;
$\beta_{k}$ is a constant weight coefficient satisfies $0 < \uline{\beta} \leq \beta_{k} \leq \bar{\beta} < \infty$.
\end{algorithmic}
\end{breakablealgorithm}

\begin{note}\label{remark}
Similar to the classical WLS algorithm,  $P_{k}$ can be expressed as $P_{k} = (\sum_{l=1}^{k} \beta_{l}\phi_{l}\phi_{l}^{T} + P_{0}^{-1})^{-1}$ (Eq. (6) in \cite{Guo123}) and $ \hat{\gamma}_{k}$ could be represented as $\hat{\gamma}_{k} = P_{k}( \sum_{l=1}^{k} \beta_{l}s_{l}\phi_{l}) + P_{k}P_{0}^{-1}\hat{\gamma}_{0}$ (Eq. (11.19) in \cite{ljung1987theory}), which indicates that $\hat{\gamma}_{k}$ given by the WLS-type algorithm is actually the minima of $J_{k}+ (\hat{\gamma}_{k}- \hat{\gamma}_{0})^{T}P_{0}^{-1}(\hat{\gamma}_{k} - \hat{\gamma}_{0})$, where the term $(\gamma - \hat{\gamma}_{0})^{T}P_{0}^{-1}(\gamma - \hat{\gamma}_{0})$ only depends on the initial values and does not affect the asymptotic properties.
\end{note}

The WLS-type algorithm has the following properties.

\begin{theorem}\label{mai_the_e}
Under Assumptions \ref{ass_a} and \ref{ass_b}, $\hat{\theta}_{k}$ given by the WLS-type algorithm has the following properties:

\romannumeral1) Almost sure convergence: $\hat{\theta}_{k}$ converges to $\theta$ in the almost sure sense with a convergence rate of $O(\sqrt{\log \log k/k})$, i.e.,
\begin{align}
\label{theta6543}\left\Vert \tilde{\theta}_{k} \right\Vert = O \left( \sqrt{\frac{\log\log k}{k}} \right), \quad \text{a.s.},
\end{align}
where $\tilde{\theta}_{k} \triangleq  \hat{\theta}_{k} - \theta$ is the estimation error of $\theta$.

\romannumeral2) $L^{p}$ convergence: $\hat{\theta}_{k}$ converges to $\theta$ in the $L^{p}$ sense with a convergence rate of $O(1/k^{p/2})$, i.e., 
\begin{align}
\label{the_1_21}& \mathbb{E}\left[\left\Vert\tilde{\theta}_{k}\right\Vert^{p} \right] = O\left(\frac{1}{k^{p/2}}\right),
\end{align}
where $p$ is an arbitrary positive integer.
\end{theorem}

The proof of Theorem \ref{mai_the_e} is supplied in Section \ref{sectiongamma}.

\begin{note}
This algorithm is primarily applicable in both noisy and noise-free cases, with the key technical aspect being that the probabilistic relationship (\ref{pri_2}) holds in both cases. 
As a result, the WLS-type criterion (\ref{ujijo2}) naturally applies to both cases. 
Furthermore, although $\hat{\theta}_{k}$ obtained based on the WLS-type criterion is related to the unknown parameter $\rho(\delta_{y})$, the ML-type algorithm for estimating $\delta_{y}$ can be conducted independently of system identification in both scenarios.
\end{note}

\begin{note}
The proposed identification algorithm provides the asymptotic optimal estimate of \( \delta_y \) under the likelihood function criterion and the asymptotic optimal solution of \( \gamma \) under the WLS-type criterion. 
However, the \( \theta \) obtained by combining these two estimates is probably not the asymptotic optimal solution of either the likelihood function criterion or the WLS-type criterion. 
Instead, it offers a feasible solution.
\end{note}

\section{Extension to dynamic systems} \label{sec_f}

The preceding theory has been developed for static linear systems, in which the output depends solely on the current input and is independent of past inputs. This section extends the proposed identification framework to handle high-dimensional OE systems, addressing both noisy and noise-free cases.

\subsection{Problem Formulation}\label{SectionVA}

Consider a high-dimensional OE system described by
\begin{align}\label{uiuiu}
y_k = \phi_{k}^{T}B(q)/A(q) + d_k,\quad k \geq 1,
\end{align}
where $k$ is the time index; $\phi_{k} \in \mathbb{R}^{n}$ is the $n$-dimensional system input; $d_{k} \in \mathbb{R}$ is the potential system noise; $A(q)$ is the $n_a$-th order scalar-valued polynomial and $B(q)$ is the $n_b$-th order $n$-dimensional vector-valued polynomial, both expressed in terms of unit backward shift operator $q^{-1}: q^{-1}y_{k} = y_{k-1}$ as 
$ A(q) = 1 + a_1 q^{-1} + \ldots + a_{n_a} q^{-n_a}$ and  $B(q) = b_0 + b_1 q^{-1} + \ldots + b_{n_b} q^{-n_b},$ $a_{1}, a_{2}, \ldots, a_{n_a} \in \mathbb{R}$ and $b_{0}, b_{1}, \ldots, b_{n_b} \in \mathbb{R}^{n}$ are unknown parameters; 
the system output $y_{k}$ cannot be exactly measured and can only be measured by quantized observations $s_{k} = Q(y_{k})$ with the quantizer $Q(\cdot)$ defined by (\ref{model_b}). 
We stipulate that  $\phi_{k} = \mathbf{0}_{n}$ for $k \leq 0$; \( A(q) \) has no poles on or outside the unit circle; $A(q)$ and $B(q)$ are coprime; $a_{n_a} \neq 0$; $b_{n_b} \neq \mathbf{0}_{n}$.

The goal of this section is to estimate unknown parameter $\theta^{*} = [a_{1},a_{2},\ldots,a_{n_a},b_{0}^{T},b_{1}^{T},\ldots,b_{n_b}^{T}]^{T} \in \mathbb{R}^{n_a+n(n_b+1)}$ based on the input $\{\phi_{k}\}_{k=1}^{\infty}$ and the quantized observation $\{ s_{k} \}_{k=1}^{\infty}$ in both noisy and noise-free cases.

\begin{note}
The OE system model generalizes the finite impulse response (FIR) system model commonly studied in quantized identification \cite{mestrah2023identification}, as the FIR system corresponds to the special case of the OE model with $A(q) = 1$.
\end{note}

\subsection{Algorithm design and its properties}

To identify the OE system, we employ an approach inspired by Durbin's method (DM) in \cite{durbin1960fitting}, which transforms the OE system identification problem into a dynamic regression system identification problem.
By incorporating a feedback structure through $A(q)$, the OE model (\ref{uiuiu}) exhibits infinite impulse response (IIR) dynamics with long memory behavior:
\begin{align}\label{IIRsystem}
y_k = \phi_k^T H(q) + d_k, \quad k \geq 1,
\end{align}
where $H(q) = B(q)/A(q) = \sum_{i=0}^{\infty} h_i q^{-i}$, and $h_0, h_1, \ldots \in \mathbb{R}^n$ are the impulse response vectors.
Under Assumptions \ref{ass_a} and \ref{ass_b}, it follows that
\begin{align*}
\mathbb{E}\left[ \phi_{k-i}y_{k} \right] = Hh_{i}, \quad i = 0,1,2,\ldots,k-1.
\end{align*}
Note that the OE system (\ref{uiuiu}) can also be described by
\begin{align}\label{uiuiu2}
A(q)y_k = \phi_{k}^{T}B(q)+ A(q)d_k,
\end{align}
and $H > 0$.
Then, multipling both sides of (\ref{uiuiu2}) by \(\phi_{k-i}\) and taking the expectation, we can obtain
\begin{align}
\label{hh1} & h_{i} + a_{1}h_{i-1} + \ldots + a_{n_a}h_{i-n_a} = b_{i}, \,  i = 0, 1,\ldots,  n_b, \\
& h_{i} + a_{1}h_{i-1} + \ldots + a_{n_a}h_{i-n_a} = \mathbf{0}_{n}, \,  i \geq n_b +1, \nonumber
\end{align}
where we stipulate that $h_{i} = 0$ for $i < 0$.
It follows that
\begin{align*}
\Gamma [a_{1},a_{2},\ldots,a_{n_a}]^{T} = [-h_{n_b+1}^{T},-h_{n_b+2}^{T},\ldots,-h_{\kappa}^{T}]^{T},
\end{align*}
where the positive integer $\kappa \geq n_a + n_b$ and 

$$
\Gamma =
 \begin{bmatrix}
h_{n_b} & h_{n_b-1} &\cdots  &h_{n_b+1- n_a}\\
h_{n_b+1}&h_{n_b}&\cdots  &h_{n_b+2-n_a} \\
\vdots & \vdots &  &  \vdots \\
h_{\kappa-1}& h_{\kappa-2}&\cdots  &h_{\kappa- n_a}
\end{bmatrix} \in \mathbb{R}^{n(\kappa-n_b)\times n_a}.
$$

\begin{proposition} \label{pri_the_uu} 
The rank of $\Gamma$ is equal to $n_a$.
\end{proposition}

The proof of Proposition \ref{pri_the_uu} is supplied in Section \ref{app_2uu}.

Hence, there exists a matrix $L \in \mathbb{R}^{ n_a \times n(\kappa-n_b)}$ such that $L\Gamma = I_{ n_a}$. This relationship implies
\begin{align}\label{incur22}
\left[a_{1},a_{2},\ldots,a_{n_a}\right]^{T} = -L\left[h_{ n_b+1}^{T},h_{ n_b+2}^{T},\ldots,h_{\kappa}^{T}\right]^{T}.
\end{align}

Thus, estimating the parameter vector $\theta^*$ reduces to estimating the impulse response vector $h = [h_0^T, h_1^T, \ldots, h_{\kappa}^T]^T$. To do this, we rewrite the OE system (\ref{uiuiu}) as a high-dimensional dynamic regression system model:
\begin{align}\label{yyyy}
y_k = \varphi_k^T h + e_k, \quad k \geq 1,
\end{align}
where $\varphi_k = [\phi_k^T, \phi_{k-1}^T, \ldots, \phi_{k-\kappa}^T]^T \in \mathbb{R}^{(\kappa + 1)n}$ is the system regressor which can be regarded as the system input, and $e_k = \phi_k^T H( q) - \varphi_k^T h + d_k$ is the modeling error which can be viewed as the system noise.
Then, based on the relationship between $h$ and $\theta^*$ given in equations (\ref{hh1}) and (\ref{incur22}), the following DM-type algorithm for estimating $\theta^*$ is outlined in Algorithm \ref{apa1}.

\begin{breakablealgorithm}
\begin{algorithmic}\label{apa1}
\caption{The DM-type algorithm for estimating $\theta^{*}$}
The algorithm for estimating $\theta^{*}$ is recursively defined at any $k \geq 1$ as follows:

\textbf{Step 1: Estimate the impulse response vector $h$:} 

Apply the WLS-type algorithm to the dynamic regression system model (\ref{yyyy}) to recursively estimate $h$, yielding $\hat{h}_{k}=[(\hat{h}_{k}^{0})^{T},(\hat{h}_{k}^{1})^{T},\ldots,(\hat{h}_{k}^{\kappa})^{T}]^{T}$.

\textbf{Step 2: Estimate the original parameter $\theta^{*}$:}
\begin{align*}
& \hat{\theta}_{k}^{*} = \left[\hat{a}_{k}^{1},\hat{a}_{k}^{2},\ldots,\hat{a}_{k}^{n_a},(\hat{b}_{k}^{0})^{T},(\hat{b}_{k}^{1})^{T},\ldots,(\hat{b}_{k}^{n_b})^{T}\right], \\
& \left[\hat{a}_{k}^{1},\hat{a}_{k}^{2},\ldots,\hat{a}_{k}^{n_a}\right]^{T} = -L_{k}\left[ (\hat{h}_{k}^{n_b+1})^{T},\ldots,(\hat{h}_{k}^{\kappa})^{T}\right]^{T},  \\
& \hat{b}_{k}^{j} = \hat{h}_{k}^{j} + \sum_{i=1}^{n_a}\hat{a}_{k}^{i}\hat{h}_{k}^{j-i}, \quad j=0,1,2,\ldots,n_b, 
\end{align*}
where $\hat{\theta}_{k}^{*}$ is the estimate of $\theta^{*}$ at time $k$; $L_{k} \in \mathbb{R}^{n_a \times \kappa}$ satifies $L_{k} \hat{\Gamma}_{k}  = I_{n_a}I_{\{ \text{rank}(\hat{\Gamma}_{k}) = n_a\}}$;
$$
\hat{\Gamma}_{k} =
\begin{bmatrix}
\hat{h}^{n_b}_{k} & \hat{h}^{ n_b-1}_{k} & \cdots  & \hat{h}^{ n_b + 1 - n_a}_{k} \\
\hat{h}^{ n_b+1}_{k}& \hat{h}^{n_b}_{k} & \cdots  & \hat{h}^{ n_b + 2 - n_a}_{k} \\
\vdots & \vdots &  & \vdots \\
\hat{h}^{\kappa-1}_{k} & \hat{h}^{\kappa-2}_{k} & \cdots & \hat{h}^{\kappa-n_a}_{k}
\end{bmatrix}.
$$
\end{algorithmic}
\end{breakablealgorithm}

\begin{note}
Unlike the static linear system discussed in Section \ref{sec_b}, the system regressor/input sequence $\{\varphi_k\}_{k=1}^{\infty}$ and  system error/noise sequence $\{e_k\}_{k=1}^{\infty}$ in (\ref{yyyy}) are no longer independent, which implies that Theorems \ref{mai_the_b} and \ref{mai_the_e} cannot be directly applied.
To address this challenge, we develop a novel analytical technique by constructing strictly stationary and geometrically strongly mixing auxiliary processes, along with carefully designed error terms that remain negligible relative to the convergence rates of the auxiliary processes themselves, to characterize the original non-iid sequences.
By proving convergence properties of these auxiliary processes, this technique ultimately establishes the almost sure convergence of the DM-type algorithm, as detailed in Section \ref{app_2}.
\end{note}

The DM-type algorithm has the following properties:
\begin{theorem}\label{thea2}
Under Assumptions \ref{ass_a} and \ref{ass_b},  if we set the weight coefficient $\beta_{k}$ to be a constant, then, $\hat{\theta}_{k}^{*}$ converges to $\theta^{*}$ almost surely with a convergence rate of $O(\sqrt{\log\log k/k})$, i.e., 
\begin{align}\label{aa_pp_2}
\left\Vert \tilde{\theta}_{k}^{*} \right\Vert & = O \left( \sqrt{\frac{\log\log k}{k}} \right), \quad \text{a.s.},
\end{align}
where $\tilde{\theta}_{k}^{*} \triangleq  \hat{\theta}_{k}^{*} - \theta^{*}$ is the estimation error of $\theta^{*}$.
\end{theorem}
The proof of Theorem \ref{thea2} is supplied in Section \ref{app_2}.

\begin{note}
The almost sure convergence rate of $O(\sqrt{\log \log k / k})$, achieved by Theorem \ref{thea2}, matches the classical WLS rate under accurate observations \cite{radenkovic2000almost} and improves upon existing $O(1/k^{v})$ rate for any $v \in (0, 1/2)$ in prior works \cite{song2018recursive}, \cite{huang2024identification}, \cite{wang2023identification}.
\end{note}

\begin{note}
The parameter \( \kappa \) represents a tradeoff between the accuracy and the computational complexity of the DM-type algorithm. Specifically, \( \kappa \) is the order of the dynamic regression system (\ref{yyyy}), derived by truncating the IIR system (\ref{IIRsystem}). 
A larger \( \kappa \) reduces the loss of \( \theta^* \) information but increases the computational cost, particularly its multiplication complexity, which grows with the order \( \kappa^3 \).
Theorem \ref{thea2} provides a lower bound for \( \kappa \), i.e., when \( \kappa \geq {n_a + n_b} \), the DM-type algorithm can successfully identify the OE systems.
\end{note}

\section{Proof of the main results}\label{sec_c}

Some lemmas are collected and established first, which will be frequently used in convergence analysis.

\begin{lemma} \label{lemma_a1}(Lemma 1 and Remark 2 in \cite{yin2007asymptotically}) 
Denote \( F^{k}(\cdot) \) as the \( k \)-sample empirical distribution of the standard Gaussian random variable.  Then, there exists a stretched Brownian bridge process \( e(\cdot) \) such that \(\sqrt{k} (F^{k}(\cdot) - F(\cdot)) \) converges to \( e(\cdot) \) almost surely and $\mathbb{E}[e(x_{1})e(x_{2})] = \min \{ F(x_{1}), F(x_{2}) \} - F(x_{1})F(x_{2})$ for all $x_{1},x_{2} \in \mathbb{R}$.
\end{lemma}

\begin{lemma} \label{lemma_a2} (Eq. (11) in \cite{ce7d2241-272f-3b4f-8cdd-b71b76d94ad4}) Denote $\mu^{r}_{k}$ as the $r$-th central moment of the Binomial process with $k$ independent trials and success probability $q$. Then, for each positive integer $p$:
\[
\lim\limits_{k \to \infty} \frac{\mu^{2p}_{k}}{(2kq(q-1))^{p}} = \frac{1 \cdot 3 \cdot 5 \cdots (2p - 1)}{2^p} = O(1).
\]
\end{lemma}

\begin{lemma} \label{lemma_a3} (Theorem 4.5.2 in \cite{chung2000course}) If the random process $\{X_{n}\}$ converges to $X$ in distribution, and for some $p > 0$, $\sup_{n}\mathbb{E}[\Vert X_{n} \Vert^{p}] = O(1)$, then, for each $r < p$:
\[
\lim\limits_{k \to \infty} \mathbb{E}[\Vert X_{n} \Vert^{r}] = \mathbb{E}[\Vert X \Vert^{r}] = O(1).
\]
\end{lemma}

\begin{lemma} \label{lemma_b1} (Proposition 6 in \cite{guo2013recursive})
The projection operator $\Pi(\cdot)$ given by Algorithm 1 follows $\Vert \Pi(x_{1}) - \Pi(x_{2}) \Vert \leq \Vert x_{1} - x_{2} \Vert$ for all $x_{1},x_{2} \in \mathbb{R}^{n}$.
\end{lemma}

\begin{lemma} \label{lemma_c1} (Lemma 3.1 in \cite{latif1979moments}) If the random symmetric positive definite matrix $A \in \mathbb{R}^{n \times n}$ has the Wishart distribution as $A \sim W_{n}(k, \Sigma)$,
where $k$ represents the degrees of freedom and $\Sigma$ is the scale matrix. 
Then, for each positive integer $p$:
\[
\mathbb{E}[\Vert A^{-1}\Vert^{p}] = \frac{\Vert \Sigma^{-1} \Vert^{p} \Gamma_{n}(k/2-p)}{2^{np}\Gamma_{n}(k/2)}, \quad \frac{k}{2} > p,
\]
where $ \Gamma_{n}(x) = \pi^{n(n-1)/4} \prod_{i=1}^{n}\Gamma(x-(i-1)/2)$ for all $x > (p-1)/2$ and $\Gamma(\cdot)$ is the gamma function.
\end{lemma}

\begin{lemma} \label{lemma_d1} (Lemma 3.1 in \cite{ledoux2013probability}) 
Let \( X \) be a Gaussian variable. Then, $\mathbb{P}(\Vert X \Vert > t) = O(\exp(-t^2/(8M^2))),$
where $M$ is the median of $\|X\|$, i.e., $\mathbb{P}(\|X\| \leq M) \geq 1/2$ and $\mathbb{P}(\|X\| \geq M) \geq 1/2.$
\end{lemma}

\begin{lemma} \label{lemma_f3} (Corollary 1.1 and Theorem 6.4 in \cite{rio2017asymptotic}) 
Let $\{X_k\}_{k = -\infty}^{\infty}$ be a strictly stationary sequence of real-valued and centered random variables satisfying 
\[
\int_0^1 \alpha^{-1}(u) \mathbb{Q}^2(u) du < \infty,
\]
where $\mathbb{Q}(u) = \inf \{ \tau \geq 0 : \mathbb{P}(\vert X_{0}\vert > \tau) \leq u \}$, \begin{align}\label{alphamix}
\alpha(\tau) =& \sup_{t \in \mathbb{Z}} \sup_{A \in \mathcal{F}_{-\infty}^t, B \in \mathcal{F}_{t+\tau}^\infty} \left| \mathbb{P}(A \cap B) - \mathbb{P}(A)\mathbb{P}(B) \right|
\end{align}
is the $\alpha$-mixing coefficient of $\{X_k\}_{k = -\infty}^{\infty}$ with the $\sigma$-algebras $\mathcal{F}_{-\infty}^t = \sigma\{X_{k}: k \leq t\}$ and $\mathcal{F}_{t+\tau}^{\infty} = \sigma\{X_{k}: k \geq t+\tau\}$, and $\alpha^{-1}(u) = \inf \{ \tau \geq 0 : \alpha(\tau)   \leq u \}$ is the càdlàg inverse function of the $\alpha$-mixing coefficient $\alpha(\tau)$. Then, 
\[
\limsup_{k \to \infty} \frac{\vert \sum_{l=1}^{k} X_{l} \vert}{\sqrt{k \log \log k}} \leq 16 \sqrt{\int_0^1 \alpha^{-1}(u) \mathbb{Q}^2(u) du}, \quad \text{a.s.}
\]
\end{lemma}

\subsection{Proof of Proposition \ref{pri_the_c}}\label{j8k}

Based on the linear system (\ref{model_a})-(\ref{model_b}),  for all $ b > 0$ and $k \geq 1$, we construct an another system:
\begin{align}
\label{model_c} \bar{y}_{k} = \phi_{k}^{T}\eta +\omega_{k}, \quad \bar{s}_{k} = I_{\{ \bar{y}_{k} > 0 \}}, \quad k \geq 1,
\end{align}
where $\eta = b\theta$ and $\omega_{k} = b d_{k}$.

Then we will show the unidentifiability between the linear system (\ref{model_a})-(\ref{model_b}) and the linear system (\ref{model_c}).

$\mathbf{Case}$ $\mathbf{1:}$ The noise-free case, i.e., $\delta_{d}^{2} = 0$.

Note that $\mathbb{P}(\bar{s}_{k}=0|\phi_{k},\eta) = I_{\{\phi_{k}^{T} \eta \leq 0\}} = I_{\{\phi_{k}^{T} \theta \leq 0\}} = \mathbb{P}(s_{k}=0|\phi_{k},\theta).$ Thus, under the same arbitrary inputs, the observations of the linear system (\ref{model_a})-(\ref{model_b}) and the linear system (\ref{model_c}) are the same in the probabilistic sense, which means the linear system (\ref{model_a})-(\ref{model_b}) is unidentifiable in this case.

$\mathbf{Case}$ $\mathbf{2:}$ The noisy case, i.e., $\delta_{d}^{2} > 0$.

Note that $\mathbb{P}(\bar{s}_{k} = 0| \phi_{k} ,  \eta)= \mathbb{P}(\phi_{k}^{T} \eta + \omega_{k} \leq 0) =  F(-\phi_{k}^{T} \eta/ b \delta_{d}) =  F(-\phi_{k}^{T} \theta/\delta_{d}) = \mathbb{P}(\phi_{k}^{T} \theta + d_{k} \leq 0) =  \mathbb{P}(s_{k} = 0| \phi_{k} ,  \theta)$. The linear system (\ref{model_a})-(\ref{model_b}) is unidentifiable in this case.

In summary, the linear system (\ref{model_a})-(\ref{model_b}) is unidentifiable with $m=1$ and $C_{1}=0$.

The proof is completed.

\subsection{Proof of Proposition \ref{pri_the_b}}\label{j7k}

From $\mathbb{E}[\phi_{k}y_{k}] = \mathbb{E}[\phi_{k}\phi_{k}^{T}\theta + \phi_{k}d_{k}] = H\theta$, we have $\mathbb{E}[ (\phi_{k} - H\theta y_{k}/\delta_{y}^{2}) y_{k}] = 0$.
Noticing both $\phi_{k} - H\theta y_{k}/\delta_{y}^{2}$ and $y_{k}$ are Gaussian, we conclude that $\phi_{k} - H\theta y_{k}/\delta_{y}^{2}$ and $y_{k}$ are independent. Thus we have
\begin{align}\label{equ_a}
\mathbb{E}\left[ \phi_{k} | y_{k}\right] =  \frac{H\theta y_{k}}{\delta_{y}^{2}}.
\end{align}
Besides, from (\ref{model_b}), it follows that
\begin{align}\label{equ_b}
\mathbb{E}[ s_{k}y_{k}] & = \sum\limits_{i=0}^{m} i \int_{C_{i}}^{C_{i+1}} \frac{x}{\sqrt{2\pi}\delta_{y}}\exp\left(-\frac{x^{2}}{2\delta_{y}^{2}}\right) dx \nonumber \\
& = \sum_{i=0}^{m}  \frac{i\delta_{y}^{2}}{\sqrt{2 \pi} \delta_{y}} \left( \exp\left(-\frac{C_{i}^{2}}{2\delta_{y}^{2}}\right)  - \exp\left(-\frac{C_{i+1}^{2}}{2\delta_{y}^{2}}\right)\right)  \nonumber \\
& =  \sum_{i=1}^{m}  \frac{\delta_{y}^{2}}{\sqrt{2 \pi} \delta_{y}} \exp\left(-\frac{C_{i}^{2}}{2\delta_{y}^{2}} \right) = \delta_{y}^{2} \rho\left(\delta_{y}\right).
\end{align}
Thus, by (\ref{equ_a}) and (\ref{equ_b}), we have $\mathbb{E}[ s_{k}\phi_{k}] = \mathbb{E}[ \mathbb{E}[s_{k}\phi_{k} | y_{k} ] ] = \mathbb{E}[s_{k}  \mathbb{E} [\phi_{k} | y_{k} ] ] = \mathbb{E}[ s_{k}y_{k}] H \theta/\delta_{y}^{2} = \rho(\delta_{y}) H \theta$.

The proof is completed.

\subsection{Proof of Proposition \ref{pri_the_oi}} \label{iu2}

By Lemma \ref{lemma_a1}, since the cumulative sum \( \sum_{i=0}^{j-1} S_k^i \) corresponds to \( F^k(C_j / \delta_y) \), for \( j = 1, 2, \dots, m \), one can get  
\begin{align}\label{uoy}
\lim_{k \to \infty} \sqrt{k} \left( \sum_{i=0}^{j-1} S_k^i - F_j \right) = e\left(\frac{C_j}{\delta_y}\right), \quad \text{a.s.}
\end{align} 
Note that $\mathbb{E}[I_{\{s_{l}=i\}}] = F_{i+1} - F_{i}$ and $\mathbb{E}[(I_{\{s_{l}=i\}} - (F_{i+1} - F_{i}))^{2}] \leq  1$, where $l = 1,2,\ldots ,k$ and $i = 0,1,2,\ldots ,m-1$.
Then, by the law of the iterated logarithm \cite{chung2000course}, for $i = 0,1,2,\ldots ,m-1$, we have
\[
\left\vert S_{k}^{i} - \left(F_{i+1}-F_{i}\right) \right\vert = O\left(\sqrt{\frac{\log\log k}{k}}\right), \quad \text{a.s}.
\]
Note that $C_{0} = -\infty$ and $F(C_{0}/\delta_y)=0$. It holds that
\begin{align}\label{uuu7777}
\left\vert\sum\limits_{i=0}^{j-1} S_{k}^{i} - F_j \right\vert = O\left(\sqrt{\frac{\log\log k}{k}}\right), \quad \text{a.s}.
\end{align}

Let \( g(\cdot) \) denote the derivative of \( F^{-1}(\cdot) \). By the mean value theorem \cite{zorich2016mathematical}, there exists a point \( \xi_k^j \) lying between \( \sum_{i=0}^{j-1} S_k^i \) and \( F_j \), such that $F^{-1} ( \sum_{i=0}^{j-1} S_k^i ) - C_j/\delta_y = g( \xi_k^j)(  \sum_{i=0}^{j-1} S_k^i - F_j)$.
Define \( \tilde{\delta}_k^j \triangleq \hat{\delta}_k^j - \delta_y \). 
Then, for nonzero \( C_j \), using (\ref{ytrt2}), we obtain
\begin{align}\label{ppp222}
\tilde{\delta}_k^j  
&= \frac{\delta_y }{F^{-1} (  \sum_{i=0}^{j-1} S_k^i )} \left( \frac{C_j}{\delta_y} - F^{-1} \left(  \sum_{i=0}^{j-1} S_k^i \right)  \right) \nonumber \\
&= -\frac{\delta_y g(\xi_k^j) ( \sum_{i=0}^{j-1} S_k^i - F_j )}{F^{-1} (  \sum_{i=0}^{j-1} S_k^i )}.
\end{align}
Note that $g(\cdot) = 1/f(\cdot)$ and $\xi_{k}^{j}$ is between $\sum_{i=0}^{j-1} S_{k}^{i}$ and $F_j$. Then, by (\ref{uuu7777}), one can get
\begin{align}
\label{huybi} & \lim\limits_{k \to \infty} g(\xi_{k}^{j}) = g(F_j) = \frac{1}{f(C_{j}/\delta_{y})}, \quad \text{a.s}, \\
\label{huybi2} & \lim_{k \to \infty} F^{-1} \left(  \sum_{i=0}^{j-1} S_k^i  \right) = \frac{C_j}{\delta_y}, \quad \text{a.s}.
\end{align}
Thus, by substituting (\ref{uoy}), (\ref{huybi}) and (\ref{huybi2}) into (\ref{ppp222}), we obtain
\begin{align}\label{y78y}
\lim\limits_{k \to \infty} \sqrt{k}\tilde{\delta}_{k}^{j} = -\frac{e(C_{j}/\delta_{y})}{f_{j}}, \quad \text{a.s}.
\end{align}

Besides, for each positive integer $p$, by (\ref{ppp222}) and  H\"{o}lder's inequality \cite{ash2014real},  it holds that
\begin{align}\label{ppp222333}
\mathbb{E}\left[\left\vert \sqrt{k} \tilde{\delta}_k^j \right\vert^{2p} \right] \leq & k^{p}\left\vert \delta_y \right\vert^{2p} \left\vert g(\xi_k^j) \right\vert^{2p} \sqrt{\mathbb{E}\left[ \left\vert \sum_{i=0}^{j-1} S_k^i - F_j \right\vert^{4p}  \right]} \nonumber \\
& \times \sqrt{\mathbb{E}\left[ \left\vert \frac{1}{F^{-1}  ( \sum_{i=0}^{j-1} S_k^i ) }\right\vert^{4p} \right]}
\end{align}
Since $\xi_{k}^{j}$ is between $\sum_{i=0}^{j-1} S_{k}^{i}$ and $F_j$, by the mean value theorem \cite{zorich2016mathematical}, $\xi_{k}^{j}$ is constrained such that its upper bound is less than $1$ and its lower bound is greater than $0$. Then, by $g(\cdot) = 1/f(\cdot)$, we have
\begin{align}\label{ppp2223331}
\left\vert g(\xi_k^j) \right\vert^{2p} = O(1).
\end{align}
Since $\sum_{i=0}^{j-1} S_k^i = \sum_{i=1}^{k}I_{\{ s_{i} < j \}}/k$ and $\mathbb{E}[\vert\sum_{i=1}^{k}I_{\{ s_{i} < j \}}-kF_{j}\vert^{2p}]$ is the $2p$-th central moment of Binomial distribution, by Lemma \ref{lemma_a2}, one can get
\begin{align}\label{ppp2223332}
\mathbb{E}\left[ \left\vert \sum_{i=0}^{j-1} S_k^i - F_j \right\vert^{4p}  \right] = O\left(\frac{1}{k^{2p}}\right).
\end{align}
Now we examine the last term on the right side of (\ref{ppp222333}). For $l = 0,1,2,\ldots,k$, it holds that $\mathbb{P}(\sum_{i=0}^{j-1} S_{k}^{i} = l/k) =  {{k}\choose{l} } F_{j}^{l}(1-F_{j})^{k-l}$.
Without loss of any generality, we assume \( C_j < 0 \). 
Then, for a positive constant $c_{s} <  1/2 - F_{j}$, we have
\begin{align}\label{Pkui}
\mathbb{P}\left(\sum_{i=0}^{j-1} S_{k}^{i} = \frac{l}{k}\right) = & {{k}\choose{l} } \left(\frac{1}{2}-c_{s}\right)^{l}\left(\frac{1}{2}+c_{s}\right)^{k-l} \nonumber \\
& \times \left( \frac{F_{j}}{1/2-c_{s}}\right)^{l}\left( \frac{1-F_{j}}{1/2+c_{s}}\right)^{k-l}.
\end{align} 
Define the function $h(x) \triangleq l\log x + (k-l)\log (1-x)$ for all $x \in \mathbb{R}$ and $\dot{h}(\cdot)$ is the derivative of \( h(\cdot) \).
Then, using the mean value theorem \cite{zorich2016mathematical}, there exists \( q_{1} \in ( F_{j}, 1/2-c_{s}) \) such that $\log (( F_{j}/(1/2-c_{s}))^{l}( (1-F_{j})/(1/2+c_{s}))^{k-l}) = h(F_{j}) - h(1/2-c_{s}) = \dot{h}(q_{1})(F_{j}-(1/2-c_{s}))$.

For all $l > k/2-c_{s}k$, we have $\dot{h}(q_{1}) = (l-kq_{1})/(q_{1}(1-q_{1})) > (1/2-c_{s}-q_{1})k/(q_{1}(1-q_{1}))$.
Define $\varrho \triangleq -(F_{j}-(1/2-c_{s}))(1/2-c_{s}-q_{1})/(q_{1}(1-q_{1}))$.
Then, one can get $\varrho > 0$ and $( F_j/(1/2-c_{s}))^l ( (1-F_j)/(1/2+c_{s}))^{k-l} = O( \exp( -\varrho k ) )$.

Note that $\binom{k}{l} (1/2-c_{s})^l (1/2+c_{s})^{k-l} \leq \sum_{l=0}^k \binom{k}{l} (1/2-c_{s})^l (1/2+c_{s})^{k-l} = 1$. Thus, by (\ref{Pkui}), for $\vert l/k -1/2 \vert < c_{s}$, it holds that
\begin{align}\label{SS11}
\mathbb{P}\left(\sum_{i=0}^{j-1} S_{k}^{i} = \frac{l}{k}\right) = O\left( \exp\left( -\varrho k \right) \right).
\end{align}

In addition, using the mean value theorem \cite{zorich2016mathematical}, there exists a point \( q_{2} \) lying between $l/k$ and $1/2$ such that $F^{-1}( l/k) -F^{-1}( 1/2) = g(q_{2})(l/k-1/2)$.
Since $g(\cdot) = 1/f(\cdot)$, we have $\vert F^{-1}(l/k) -F^{-1}(1/2) \vert \geq \sqrt{2\pi} \vert l/k -1/2\vert$.
Thus, for $l/k \neq 1/2$, it holds that $\vert F^{-1}(l/k) - F^{-1}(1/2)\vert \geq \sqrt{2\pi}/(2k)$.
Since the point $\sum_{i=0}^{j-1} S_{k}^{i} =1/2$ has been modified, by (\ref{SS11}), we have
\begin{align*}
& \mathbb{E}\left[ \left\vert \frac{1}{F^{-1}  ( \sum_{i=0}^{j-1} S_k^i ) }\right\vert^{4p} I_{ \{ \vert \sum_{i=0}^{j-1} S_k^i - 1/2\vert < c_{s}\}}\right] \\
= &  \sum_{1/2k \leq \vert l/k -1/2 \vert < c_{s}} \frac{\mathbb{P}(\sum_{i=0}^{j-1} S_{k}^{i} = l/k)}{\vert F^{-1}  (l/k)\vert^{4p}} +  \frac{\mathbb{P}(\sum_{i=0}^{j-1} S_{k}^{i} = 1/2)}{\vert F^{-1} ( c^{*})\vert^{4p}}\\
= & O(k^{4p + 1} \exp\left( -\varrho k \right) ) = O(1).
\end{align*}
Besides, for $\vert \sum_{i=0}^{j-1} S_k^i - 1/2 \vert \geq c_{s}$, $\vert F^{-1} ( \sum_{i=0}^{j-1} S_k^i )\vert^{-4p}$ is constrained such that its lower bound is greater than $0$, which indicates that $\mathbb{E}[ \vert F^{-1} ( \sum_{i=0}^{j-1} S_k^i )\vert^{-4p} I_{ \{ \vert \sum_{i=0}^{j-1} S_k^i - 1/2 \vert \geq c_{s}\}}] = O(1)$.
Therefore, we have
\begin{align}\label{ppp2223333}
\mathbb{E}\left[ \left\vert \frac{1}{F^{-1}  ( \sum_{i=0}^{j-1} S_k^i ) }\right\vert^{4p} \right] = O(1).
\end{align}
Thus, by substituting (\ref{ppp2223331}), (\ref{ppp2223332}) and (\ref{ppp2223333}) into (\ref{ppp222333}), we obtain
\begin{align}\label{pth}
\mathbb{E}\left[\left\vert\sqrt{k}\tilde{\delta}_{k}^{j}\right\vert^{2p}\right] = O(1).
\end{align}
Consequently, in the case of $\forall C_j \neq 0$, one can get 
\begin{align}\label{Delya}
\mathbb{E}\left[\left\Vert\sqrt{k}\tilde{\Delta}_{k}\right\Vert^{4}\right] = O(1).
\end{align}

In addition, by Lemma \ref{lemma_a1}, for $u \leq v$, we have
\begin{align}\label{yyyy9999}
\mathbb{E}\left[e\left(\frac{C_{u}}{\delta_y}\right)e\left(\frac{C_{v}}{\delta_y}\right)\right] = F_u - F_v F_v.
\end{align}
Then, by Lemma \ref{lemma_a3}, (\ref{y78y}) and (\ref{Delya}), in the case of $\forall C_j \neq 0$, one can get
\begin{align*}
\lim\limits_{k \to \infty} kV_{k} & = \lim\limits_{k \to \infty} \mathbb{E}\left[k\tilde{\Delta}_{k}\tilde{\Delta}_{k}^{T}\right] 
\\ 
& =   \left(U+G^{T}\right)^{-1}\left(W-ww^T\right)\left(U+G\right)^{-1}.
\end{align*}

In the case of $C_{j} = 0$ with $j > 1$, by the mean value theorem \cite{zorich2016mathematical}, there exists a point $\zeta_{k}^{j}$ between $F(0) - \sum_{i=1}^{j-1} S_{k}^{i}$ and $F_1$ such that 
\begin{align}\label{yuyuihi2}
\tilde{\delta}_{k}^{j} &
 =  \frac{\delta_{y}g(\zeta_{k}^{j})( F_{1} - (F(0) - \sum_{i=1}^{j-1} S_{k}^{i}))}{F^{-1}(F(0) - \sum_{i=1}^{j-1} S_{k}^{i})}.
\end{align}
Then, by $F_{j} = F(0)$ and (\ref{uoy}), it holds that $\sqrt{k}( F_{1} - (F(0) - \sum_{i=1}^{j-1} S_{k}^{i})$ converges to $e(C_{j}/\delta_{y}) - e(C_{1}/\delta_{y})$ almost surely.
Similarly with the analysis of the nonzero thresholds, for $C_{j} = 0$ with $j > 1$, one can get (\ref{pth}) holds and
\begin{align}
\label{y78y2}& \lim\limits_{k \to \infty} \sqrt{k}\tilde{\delta}_{k}^{j}  = \frac{e(C_{j}/\delta_y)-e(C_{1}/\delta_y)}{f_{1}}, \quad \text{a.s}.
\end{align}
Thus, in the case of $C_{j} = 0$ with $j > 1$, we have (\ref{Delya}) holds. 

Define $U_{j} \triangleq \text{diag}(f_{1},f_{2},\ldots ,f_{j-1},-f_{1},f_{j+1},\ldots ,f_{m}) \in \mathbb{R}^{m\times m}$ and a positive definitive matrix $Y \in \mathbb{R}^{m \times m}$ satifies 
\begin{align*}
Y_{p,q} \triangleq &  \left(W -ww^T\right)_{p,q}, \, p \neq j, q \neq j, \,1 \leq p,q \leq m; \\
Y_{p,j} = & Y_{j,p} \triangleq F_{p} - F_{j}F_{p} - (F_{1} - F_{1}F_{p}) \\
= &  \left(W -ww^T\right)_{j,q} -  \left(W -ww^T\right)_{1,p}, \, p \neq j, \, 1 \leq p \leq m;  \\
Y_{j,j} \triangleq & F_{j} - F_{j}^{2}-F_{1}+2F_{1}F_{j} - F_{1}^2  \\
= &  \left(W -ww^T\right)_{j,j} - 2 \left(W -ww^T\right)_{1,j} +  \left(W -ww^T\right)_{1,1}.
\end{align*}
Then, in the case of $C_{j} = 0, j>1$, by Lemma \ref{lemma_a3}, (\ref{y78y}), (\ref{Delya}),  (\ref{yyyy9999}) and (\ref{y78y2}), we have
\begin{align}\label{uhiu23}
\lim\limits_{k \to \infty} kV_{k} & = U_{j}^{-1}Y U_{j}^{-1}.
\end{align}

Define $Z \in \mathbb{R}^{m \times m}$ satifing 
\begin{align*}
Z_{u,v} =
\left\{
\begin{array}{ll}
1,  & \text{if} \, 1 \leq u = v \leq m \, \text{or} \, u = 1, v = j, \\
0, & \text{otherwise}.
\end{array} 
\right.
\end{align*}
Then, one can get $Z^{T}YZ = W-ww^{T}$.
Thus, $Z^{T}YZ = W-ww^{T}$. 
Besides, by $U_{j}Z = U+G$ and (\ref{uhiu23}), in the case of $C_{j} = 0, j>1$, we obtain
 \begin{align*}
\lim\limits_{k \to \infty} kV_{k}
& = \left(Z^{T}U_{j}\right)^{-1}Z^{T}YZ\left(U_{j}Z\right)^{-1} \\
& = \left(U+G^{T}\right)^{-1}\left(W -ww^T\right)(U+G)^{-1}.
\end{align*}

Similarly, (\ref{pri_3}) can also be proven in the case of $C_{1} = 0$.

The proof is completed.

\subsection{Proof of Theorem \ref{mai_the_b}}\label{sectiondeltay}

This proof will be divided into the following two parts.

$\mathbf{Part}$ $\mathbf{1:}$ Represent the CR lower bound.

By $\partial F_{i} / \partial \delta_{y}= -f_i$ and (\ref{ttr}), in the case of $\forall C_j \neq 0$, it holds that
\begin{align*}
\frac{ \partial^{2} l_{k} }{ k\partial \delta_{y}^{2} }  = 
& - \frac{S_{k}^{0} \dot{f}_{1}}{F_{1}}
- \frac{S_{k}^{0}f_{1}^{2}}{F_{1}^2} +  \frac{S_{k}^{m}\dot{f}_{m}}{1- F_{m}}  -  \frac{S_{k}^{m}f_{m}^2}{(1- F_{m})^2} \nonumber \\ & + \sum\limits_{i=1}^{m-1}  \frac{S_{k}^{i}(\dot{f}_{i} - \dot{f}_{i+1})}{F_{i+1} - F_{i}}  -  \frac{S_{k}^{i}(f_{i} - f_{i+1})^2}{(F_{i+1} - F_{i})^2} \nonumber \\
= & - \sum\limits_{i=0}^{m} S_{k}^{i} \frac{\dot{\tilde{f}}_{i}}{\tilde{F}_{i}} - S_{k}^{i} \frac{\tilde{f}_{i}^2}{\tilde{F}_{i}^2}.
\end{align*}
Similarly, in the case of $\exists C_j = 0$, one can get
\begin{align*}
\frac{ \partial^{2} l_{k} }{ k\partial \delta_{y}^{2} }  
= & -\sum\limits_{i=0}^{j-2} S_{k}^{i} \frac{\dot{\tilde{f}}_{i}}{\tilde{F}_{i}} + S_{k}^{j-1}\frac{\dot{f}_{j-1}}{\tilde{F}_{j-1}} - S_{k}^{j}\frac{\dot{f}_{j+1}}{\tilde{F}_{j}} \nonumber  \\
 &- \sum\limits_{i=j+1}^{m} S_{k}^{i} \frac{\dot{\tilde{f}}_{i}}{\tilde{F}_{i}} + \sum\limits_{i=0}^{m} S_{k}^{i} \frac{\tilde{f}_{i}^2}{\tilde{F}_{i}^2}.
\end{align*}

From $\mathbb{E}[S_{k}^{i}] = \mathbb{E}[\sum_{l=1}^{k}I_{\{s_{l}=i\}}/k] = (F_{i+1} - F_{i}) = \tilde{F}_{i}$ in both cases, where $i = 1,2,\ldots,m$, we have
\begin{align}\label{koipo}
\sigma_{\text{CR}}(k) = -\left(\mathbb{E}\left[\frac{ \partial^{2} l_{k} }{ \partial \delta_{y}^{2} }\right]\right)^{-1} =  \frac{1}{k}\left(\sum\limits_{i=0}^{m} \frac{\tilde{f}_{i}^2}{\tilde{F}_{i}}\right)^{-1}.
\end{align}

Define $\tilde{w} \triangleq [\tilde{F}_{0},\tilde{F}_{1},\ldots ,\tilde{F}_{m-1}]^{T} \in \mathbb{R}^{m}$, $\tilde{W} \triangleq \text{diag}(\tilde{F}_{0},\tilde{F}_{1},\ldots ,\tilde{F}_{m-1})  \in \mathbb{R}^{m\times m}$. Additionally, let $\tilde{a} \triangleq (1 - \tilde{w}^{T}\tilde{W}^{-1}\tilde{w})^{-1}$.
Then, we have
\begin{align}\label{Www1}
(\tilde{W}-\tilde{w}\tilde{w}^{T})^{-1} & = \tilde{W}^{-1/2}(I_{m}-\tilde{W}^{-1/2}\tilde{w}\tilde{w}^{T}\tilde{W}^{-1/2})^{-1}\tilde{W}^{-1/2}\nonumber\\
& = \tilde{W}^{-1/2}(I_{m} + \tilde{a}\tilde{W}^{-1/2}\tilde{w}\tilde{w}^{T}\tilde{W}^{-1/2})\tilde{W}^{-1/2} \nonumber\\
& = \tilde{W}^{-1} + \tilde{a}\tilde{W}^{-1}\tilde{w}\tilde{w}^{T}\tilde{W}^{-1}.
\end{align}

Define $\tilde{U} \triangleq \text{diag}(\tilde{f}_{0},\tilde{f}_{1},\ldots,\tilde{f}_{m-1}) \in \mathbb{R}^{m\times m}$ and $\tilde{X} \triangleq \tilde{U}(\tilde{W}-\tilde{w}\tilde{w}^{T})^{-1}\tilde{U}$. 
It follows that 
\begin{align}
\tilde{X} = \tilde{U}\tilde{W}^{-1}\tilde{U} + \tilde{a}\tilde{U}\tilde{W}^{-1}\tilde{w}\tilde{w}^{T}\tilde{W}^{-1}\tilde{U}.
\end{align}
Note that $\mathbf{1}_{m}^{T}\tilde{U}\tilde{W}^{-1}\tilde{U}\mathbf{1}_{m} =  \sum_{i=0}^{m-1} \tilde{f}_{i}^2/\tilde{F}_{i},$ $\tilde{a}^{-1} = 1 - \tilde{w}^{T}\tilde{W}^{-1}\tilde{w} = 1 - \sum_{i=0}^{m-1}\tilde{F}_{i} = \tilde{F}_{m},$ and $ \tilde{w}^{T}\tilde{W}^{-1}\tilde{U}\mathbf{1}_{m} =  \sum_{i=0}^{m-1} \tilde{f}_{i} = -\tilde{f}_{m}.$
Then, we have 
\begin{align}\label{uhiuh}
\mathbf{1}_{m}^{T}\tilde{X} \mathbf{1}_{m} = \sum_{i=0}^{m} \tilde{f}_{i}^2/\tilde{F}_{i}.
\end{align} 

Define  $X \triangleq (U+G)(W -ww^T)^{-1}(U+G^{T})$ and
\[
Z_{1} \triangleq
 \begin{bmatrix}
1 & -1 &\cdots  & -1 \\
0 & 1 &\cdots  & 0 \\
\vdots & \vdots & \ddots &  \vdots \\
0 & 0 &\cdots  & 1
\end{bmatrix}, 
\]
which is formed by subtracting columns $i = 2,\ldots,m$ from the first column of the identity matrix. Since $\mathbf{1}_{m}^{T}G = [0,0,\ldots,0]  \in \mathbb{R}^{m}$, one may obtain $ \mathbf{1}_{m}^{T}X \mathbf{1}_{m} = \mathbf{1}_{m}^{T}UZ_{1}(Z_{1}^{T}WZ_{1}-Z_{1}^{T}ww^{T}Z_{1})^{-1}Z_{1}^{T}U\mathbf{1}_{m},$
and 
\begin{align*}
& \mathbf{1}_{m}^{T}UZ_{1} = \left[f_1,f_2 - f_1, \ldots ,f_m-f_1\right], \\
& w^{T}Z_{1} = \left[F_1,F_2 - F_1, \ldots ,F_m-F_1\right],   \\
& Z_{1}^{T}WZ_{1} =  
\begin{bmatrix}
F_{1} & 0 &\cdots  &0\\
0 &F_{2}-F_1&\cdots  &F_{2}-F_1 \\
\vdots &\vdots & \ddots &  \vdots \\
0 & F_{2}-F_1&\cdots  &F_{m}-F_1
\end{bmatrix}.
\end{align*}

This process can be repeated. For instance, we denote
\[
Z_{2} \triangleq Z_{1} \cdot
 \begin{bmatrix}
1 & 0 & 0 & \cdots  & 0 \\
0 & 1 & -1 & \cdots  & -1 \\
0 & 0 & 1 & \cdots  & 0 \\
\vdots &\vdots &\vdots &\ddots &  \vdots \\
0 & 0 & 0 &\cdots  & 1
\end{bmatrix}, 
\]
where the later matrix is formed by subtracting columns $i = 3,4,\ldots,m$ from the second column of $I_{m}$. After $m-1$ such elementary operations, we obtain $\mathbf{1}_{m}^{T}UZ_{m-1} = \mathbf{1}_{m}^{T}\tilde{U}$ and
\begin{align}\label{www222000111}
w^{T}Z_{m-1} = \tilde{w}, \quad Z_{m-1}^{T}WZ_{m-1} =  \tilde{W}.
\end{align}
By (\ref{koipo}) and (\ref{uhiuh}), we have
\begin{align}\label{uhiuhi}
\mathbf{1}_{m}^{T}X \mathbf{1}_{m} = \mathbf{1}_{m}^{T}\tilde{X} \mathbf{1}_{m} = \frac{1}{k\sigma_{\text{CR}}(k)}.
\end{align}

$\mathbf{Part}$ $\mathbf{2:}$ Establish the convergence properties of $\hat{\delta}_{k}$.

We first prove that $\Vert\hat{\mu}_{k}\Vert = O(1)$.

Define $\tilde{w}_{k} \triangleq [\tilde{F}_{k}^{0},\tilde{F}_{k}^{1},\ldots ,\tilde{F}_{k}^{m-1}]^{T} \in \mathbb{R}^{m}$, $\tilde{W}_{k} \triangleq \text{diag}(\tilde{F}_{k}^{0},\tilde{F}_{k}^{1},\ldots ,\tilde{F}_{k}^{m-1})  \in \mathbb{R}^{m\times m}$ and $\tilde{a}_{k} \triangleq (1 - \tilde{w}_{k}^{T}\tilde{W}_{k}^{-1}\tilde{w}_{k})^{-1} \in \mathbb{R}$, where $\tilde{f}_{k}^{l} \triangleq \hat{f}_{k}^{l+1} - \hat{f}_{k}^{l}$ and $\tilde{F}_{k}^{l} \triangleq \hat{F}_{k}^{l+1} - \hat{F}_{k}^{l}$ for $l = 0,1,2,\ldots,m$.
Then, similarly with (\ref{Www1}), one can get $
(\tilde{W}_{k}-\tilde{w}_{k}\tilde{w}_{k}^{T})^{-1} = \tilde{W}_{k}^{-1} + \tilde{a}_{k}\tilde{W}_{k}^{-1}\tilde{w}_{k}\tilde{w}_{k}^{T}\tilde{W}_{k}^{-1}$.

Since $\hat{\delta}_{k-1}$ has a positive lower bound and a bounded upper bound,  $\tilde{F}_{l}^{k}$ and $\tilde{f}_{l}^{k}$ also has positive lower bounds and bounded upper bounds.
This indicates that $\Vert \tilde{W}_{k}^{-1} \Vert$, $\vert \tilde{w}_{k} \vert$ and $\Vert \hat{U}_{k}+\hat{G}_{k} \Vert$ are all bounded.
Besides, since $\tilde{a}_{k}^{-1} = \tilde{F}_{m}^{k}$, $\vert \tilde{a}_{k} \vert$ is bounded.  
Thus, from $\Vert (\tilde{W}_{k}-\tilde{w}_{k}\tilde{w}_{k}^{T})^{-1} \Vert \leq \Vert \tilde{W}_{k}^{-1} \Vert +  \Vert \tilde{a}_{k} \Vert\Vert \tilde{W}_{k}^{-1} \Vert^{2}\Vert \tilde{w}_{k} \Vert^{2}$, we know $\Vert (\tilde{W}_{k}-\tilde{w}_{k}\tilde{w}_{k}^{T})^{-1} \Vert$ is bounded.

Note that $(\hat{W}_{k}-\hat{w}_{k}\hat{w}_{k}^T)^{-1} = Z_{m-1}(\tilde{W}_{k}-\tilde{w}_{k}\tilde{w}_{k}^{T})^{-1}Z_{m-1}^{T}$. Then, we have $\Vert (\hat{W}_{k}-\hat{w}_{k}\hat{w}_{k}^T)^{-1} \Vert$ is bounded.

Consequently, $\Vert(\hat{U}_{k}+\hat{G}_{k})(\hat{W}_{k}-\hat{w}_{k}\hat{w}_{k}^T)^{-1}(\hat{U}_{k}+\hat{G}_{k}^{T}) \Vert$ is bounded.

In addition, similarly with (\ref{uhiuh}) and (\ref{uhiuhi}), we have $\mathbf{1}_{m}^{T}(\hat{U}_{k}+\hat{G}_{k})(\hat{W}_{k}-\hat{w}_{k}\hat{w}_{k}^T)^{-1}(\hat{U}_{k}+\hat{G}_{k}^{T})\mathbf{1}_{m} =  \sum_{i=0}^{m} (\tilde{f}_{k}^{i})^2/\tilde{F}_{k}^{i}$, where the positive lower bound and the bounded upper bound of $\hat{\delta}_{k-1}$ ensures that $\vert (\mathbf{1}_{m}^{T}(\hat{U}_{k}+\hat{G}_{k})(\hat{W}_{k}-\hat{w}_{k}\hat{w}_{k}^T)^{-1}(\hat{U}_{k}+\hat{G}_{k}^{T})\mathbf{1}_{m})^{-1} \vert$ is also bounded.
Therefore, by (\ref{www23451}), we have
\begin{align}\label{boundmu}
\Vert \hat{\mu}_{k} \Vert = O(1).
\end{align}

For \( C_j \neq 0 \), according to (\ref{uoy}), (\ref{uuu7777}) and (\ref{y78y}), we have
\begin{align}\label{yuyugu}
\left\vert \tilde{\delta}_{k}^{j} \right\vert = O\left( \sqrt{\frac{\log\log k}{k}} \right), \quad \text{a.s.}
\end{align}
For \( C_j = 0 \) with \( j > 1 \), a similar argument based on (\ref{yuyuihi2}) demonstrates that (\ref{yuyugu}) remains valid. Likewise, in the case where \( C_1 = 0 \), the same conclusion holds, establishing that (\ref{yuyugu}) is satisfied in all cases. 

 Since the positive constant $c$ is sufficiently small, we assume $\delta_{y} \in (c,1/c)$. Then, by (\ref{www23456}) and Lemma \ref{lemma_b1}, we have $\vert\tilde{\delta}_{k}\vert \leq \Vert\hat{\mu}_{k} [\tilde{\delta}_{k}^{1}, \tilde{\delta}_{k}^{2}, \ldots ,\tilde{\delta}_{k}^{m}]\Vert \leq \Vert\hat{\mu}_{k} \Vert \Vert[\tilde{\delta}_{k}^{1}, \tilde{\delta}_{k}^{2}, \ldots ,\tilde{\delta}_{k}^{m}]\Vert$.
By (\ref{yuyugu}), one can get (\ref{y7tgyghuh}) holds, i.e., the almost sure convergence rate is $O(\sqrt{\log\log k/k})$.

Define $\acute{\delta}_{k} \triangleq [\hat{\delta}_{k}^{1}, \hat{\delta}_{k}^{2}, \ldots ,\hat{\delta}_{k}^{m}]\hat{\mu}_{k}$ and $\bar{\delta}_{k} \triangleq \acute{\delta}_{k} - \delta_{y}$.
Note that for \( C_j \neq 0 \) and \( C_j = 0 \) with \( j > 1 \), (\ref{pth}) remains valid. 
Similarly, the same conclusion holds where \( C_1 = 0 \).  
By applying (\ref{pth}), (\ref{boundmu}) and the Rosenthal inequality \cite{rosenthal1970subspaces}, we obtain $\mathbb{E}[\vert \bar{\delta}_{k} \vert^{2p} ] = O( \sum_{j=1}^{m} \mathbb{E}[\vert \tilde{\delta}_{k}^{j} \vert^{2p}] + (\sum_{j=1}^{m} \mathbb{E}[\vert \tilde{\delta}_{k}^{j} \vert^{2}] )^p ) = O(1/k^{p})$.
Then, by Lemma \ref{lemma_b1} and H\"{o}lder's inequality \cite{ash2014real}, we have $\mathbb{E}[\vert \tilde{\delta}_{k} \vert^{p} ]  \leq \mathbb{E}[ \vert \bar{\delta}_{k} \vert^{p}] \leq (\mathbb{E}[\vert \bar{\delta}_{k}\vert^{2p}])^{1/2} = O(1/k^{p/2})$,
which directly indicates that (\ref{the_1_2}) holds, i.e., the $L^{p}$ convergence rate is $O(k^{p/2})$.

Define $\uline{c} \triangleq \min\{ \delta_{y} -c, 1/c -\delta_{y} \}$.
Then, by the Markov inequality \cite{P}, one can get $\mathbb{P}( \acute{\delta}_{k} \notin [c,1/c] ) \leq \mathbb{P}(\vert \bar{\delta}_{k}\vert > \uline{c} ) \leq \mathbb{E}[\vert \bar{\delta}_{k} \vert^{2}]/\uline{c}^{2} = O(1/k).$
Thus, by H\"{o}lder's inequality \cite{ash2014real}, we have
$k\mathbb{E}[ \vert \bar{\delta}_{k} - \tilde{\delta}_{k} \vert^{2} ] \leq k (\mathbb{E}[ \vert \bar{\delta}_{k} - \tilde{\delta}_{k} \vert^{4} ] \mathbb{P}( \acute{\delta}_{k} \notin [c,1/c] ))^{1/2} \leq k (2(\mathbb{E}[ \vert \bar{\delta}_{k}\vert^{4}]+  \mathbb{E}[\tilde{\delta}_{k} \vert^{4} ] )\mathbb{P}( \acute{\delta}_{k} \notin [c,1/c] ))^{1/2} =  O(1/k^{1/2})$.
Thus, one can get
\begin{align}\label{uu77yy22}
& \lim\limits_{k \to \infty} \mathbb{E}\left[\left\vert \sqrt{k}\tilde{\delta}_{k}\right\vert^{2} \right] - \mathbb{E}\left[\left\vert \sqrt{k}\bar{\delta}_{k}\right\vert^{2} \right]\nonumber \\
= & \lim\limits_{k \to \infty}  2k\mathbb{E}\left[\bar{\delta}_{k}  \left(\tilde{\delta}_{k} - \bar{\delta}_{k}  \right) \right]  + k\mathbb{E}\left[ \left( \bar{\delta}_{k} - \tilde{\delta}_{k} \right)^{2} \right]  
=  0.
\end{align}
Besides, by (\ref{www23451}) and (\ref{y7tgyghuh}), we obtain $(\hat{U}_{k}+\hat{G}_{k})(\hat{W}_{k}-\hat{w}_{k}\hat{w}_{k}^T)^{-1}(\hat{U}_{k}+\hat{G}_{k}^{T})$ converges to $X$ almost surely.
Then, by (\ref{www23451}), we have
\[
\lim\limits_{k \to \infty} \left\vert \sqrt{k}\bar{\delta}_{k} \right\vert^2 - k\frac{\mathbf{1}_{m}^{T}X \tilde{\Delta}_{k}\tilde{\Delta}_{k}^{T}X\mathbf{1}_{m}}{(\mathbf{1}_{m}^{T}X\mathbf{1}_{m})^{2}} = 0, \quad  \text{a.s.}
\]
By Lemma \ref{lemma_a3} and Proposition \ref{pri_the_oi}, it holds that
\begin{align}\label{uuhh2233}
\lim\limits_{k \to \infty} \mathbb{E}\left[\left\vert \sqrt{k}\bar{\delta}_{k}\right\vert^{2} \right]-k\sigma_{\text{CR}}(k) = 0.
\end{align}
Then, by (\ref{uu77yy22}), one can get (\ref{theorem_3}) holds, i.e., the asymptotic variance is established.

The proof is completed.

\subsection{Proof of Theorem \ref{mai_the_e}}\label{sectiongamma}

This proof will be divided into the following two parts.

$\mathbf{Part}$ $\mathbf{1:}$ Establish the almost sure convergence of $\hat{\theta}_{k}$.

Define $\tilde{\gamma}_{k} \triangleq \hat{\gamma}_{k} - \gamma$ for $k \geq 0$. 
From Remark \ref{remark}, we have $\tilde{\gamma}_{k} = P_{k}(( \sum_{l=1}^{k} \beta_{l}s_{l}\phi_{l}) + P_{0}^{-1}\hat{\gamma}_{0} - P_{k}^{-1}\gamma)$ and $P_{k}=(\sum_{l=1}^{k} \beta_{l}\phi_{l}\phi_{l}^{T} + P_{0}^{-1})^{-1}$, which follows that
\begin{align}\label{iuy}
\tilde{\gamma}_{k} = P_{k}\sum\limits_{l=1}^{k} \beta_{l}\phi_{l}\left(s_{l}-\phi_{l}^{T}\gamma\right) + P_{k}P_{0}^{-1}\left(\hat{\gamma}_{0}-\gamma\right).
\end{align} 

For the simplicity of description, we denote $x_{k} \triangleq s_{k}\phi_{k}-\phi_{k}\phi_{k}^{T}\gamma$. 
From Proposition \ref{pri_the_b}, we obtain
\begin{align}\label{equ_ss}
\mathbb{E}\left[\beta_{k}x_{k}\right] =  \beta_{k}H\gamma - \beta_{k}H\gamma = \mathbf{0}_{n}.
\end{align}
Since $\{\phi_{k}\}_{k=1}^{\infty}$ is Gaussian iid, $\mathbb{E}[\Vert \phi_{l} \Vert^{4}] < \infty$. Together with $\Vert \gamma \Vert =  \Vert \sum_{i=1}^{m}  \exp(-C_{i}^{2}/(2\delta_{y}^{2})) \theta /(\sqrt{2 \pi} \delta_{y}) \Vert < \infty$, by the $C_{r}$-inequality \cite{P}, it holds that
\begin{align}\label{sdssda}
 \mathbb{E}\left[ \left\Vert \beta_{k} x_{k} \right\Vert^{4p} \right]  \leq & 2^{4p-1} \bar{\beta}^{4p} \mathbb{E}\left[ \left(s_{k}^{4p}+\left(\phi_{k}^{T}\gamma\right)^{4p}\right)\left\Vert \phi_{k} \right\Vert^{4p}\right] \nonumber \\ \leq & 2^{4p-1}\bar{\beta}^{4p} \mathbb{E}\left[m^{4p}\left\Vert \phi_{k} \right\Vert^{4p} +\left\Vert \gamma\right \Vert^{4p} \left\Vert \phi_{k}\right \Vert^{8p}\right]\nonumber\\ = & O(1).
\end{align}
Then, by the law of the iterated logarithm \cite{chung2000course}, we have 
\begin{align} \label{equ_sss}
\left\Vert \sum_{l=1}^{k}\beta_l x_{l} \right\Vert = O\left(\sqrt{k \log \log k}\right), \quad \text{a.s.}
\end{align}
By the law of the iterated logarithm \cite{chung2000course}, it holds that $\sum_{l=1}^{k}\phi_{l}\phi_{l}^{T}/k$ converges to $H$ almost surely, which indicates that $\Vert k (\sum_{l=1}^{k}\phi_{l}\phi_{l}^{T})^{-1}  \Vert$ converges to $\Vert H^{-1}\Vert$ almost surely.
Since $\Vert P_{k} \Vert = \Vert (\sum_{l=1}^{k} \beta_{l}\phi_{l}\phi_{l}^{T} + P_{0}^{-1})^{-1} \Vert \leq \Vert (\sum_{l=1}^{k}\phi_{l}\phi_{l}^{T})^{-1} \Vert/\uline{\beta} $, it holds that
\begin{align}\label{ppppi}
 \left \Vert P_{k} \right \Vert   = O\left(\frac{1}{k}\right), \quad \text{a.s.} 
\end{align}
Therefore, by substituting (\ref{equ_sss}) and (\ref{ppppi}) into (\ref{iuy}), we obtain
\begin{align}\label{m76}
\left\Vert \tilde{\gamma}_{k} \right \Vert = O\left(\sqrt{\frac{\log\log k}{k}}\right), \quad  \text{a.s.}
\end{align}
In addition, by the mean value theorem \cite{zorich2016mathematical}, there exists a point $\eta_{k}$ between $\hat{\delta}_{k}$ and $\delta_{y}$ such that $\rho(\hat{\delta}_{k}) - \rho(\delta_{y}) = \dot{\rho}(\eta_{k}) \tilde{\delta}_{k}$, where $\dot{\rho}(\cdot)$ is the derivative of $\rho(\cdot)$.
Note that $\dot{\rho}(\eta_{k}) = \sum_{i=1}^{m} ((C_{i}^{2}-\eta_{k}^{2}) \exp(- C_{i}^{2}/(2\eta_{k}^{2})))/\sqrt{2\pi}\eta_{k}^{4})$.
Since $\vert\hat{\delta}_{k}\vert$ is bounded, $\vert\dot{\rho}(\eta_{k}) \vert$ is also bounded. By Theorem \ref{mai_the_e}, we have
\begin{align}\label{uyui1}
\left\vert \rho(\hat{\delta}_{k}) - \rho(\delta_{y}) \right\vert =  O\left(\sqrt{\frac{\log \log k}{k}}\right), \quad \text{a.s.}
\end{align}

By $\tilde{\theta}_{k} = \hat{\gamma}_{k}/\rho(\hat{\delta}_{k}) - \gamma/\rho(\delta_{y})$, one can get
\begin{align}\label{jioj100}
\left\Vert \tilde{\theta}_{k} \right\Vert 
& = \left\Vert \frac{\hat{\gamma}_{k}\rho(\delta_{y}) -  \gamma\rho(\delta_{y}) + \gamma\rho(\hat{\delta}_{k}) - \gamma\rho(\hat{\delta}_{k})}{\rho(\hat{\delta}_{k})\rho(\delta_{y})} \right\Vert \nonumber \\
& \leq \left\Vert \frac{\tilde{\gamma}_{k}}{\rho(\hat{\delta}_{k})} \right\Vert + \left\Vert \frac{(\rho(\hat{\delta}_{k}) - \rho(\delta_{y}))\gamma}{\rho(\hat{\delta}_{k})\rho(\delta_{y})}  \right\Vert. 
\end{align}
Since $\hat{\delta}_{k}$ has a positive lower bound, $\rho(\hat{\delta}_{k})$ has a positive lower bound. By substituting (\ref{m76}) and (\ref{uyui1}) into (\ref{jioj100}), we obtain (\ref{theta6543}), i.e., the almost sure convergence rate is $O(\sqrt{\log\log k/k})$.

$\mathbf{Part}$ $\mathbf{2:}$ Establish the $L^{p}$ convergence of $\hat{\theta}_{k}$.

By (\ref{iuy}) and the $C_{r}$-inequality \cite{P}, we have 
\begin{align*}
\mathbb{E}\left[ \left\Vert \sqrt{k}\tilde{\gamma}_{k}\right\Vert^{2p} \right] \leq & 2^{2p-1}k^{p}\mathbb{E}\left[ \left\Vert P_{k} \right\Vert^{2p} \left\Vert P_{0} \right\Vert^{2p} \left\Vert\hat{\gamma}_{0}-\gamma \right\Vert^{2p}\right] \nonumber \\ & + 2^{2p-1} k^{p}\mathbb{E}\left[ \left\Vert P_{k} \right\Vert^{2p} \left\Vert \sum_{l=1}^{k}\beta_{l}x_{l} \right\Vert^{2p} \right].
\end{align*}
Since $\{\beta_{k}x_{k}\}_{k=1}^{\infty}$ is an independent sequence with zero mean and finite moments of arbitrary order, by the Rosenthal inequality \cite{rosenthal1970subspaces}, it holds that $\mathbb{E}[ \Vert \sum_{l=1}^{k}\beta_{l}x_{l} \Vert^{4p}] = O(\sum_{l=1}^{k}\mathbb{E}[ \Vert \beta_{l}x_{l}\Vert^{4p}] + (\sum_{l=1}^{k}\mathbb{E}[ \Vert \beta_{l}x_{l} \Vert^{2}])^{2p}) = O( k^{2p} )$.

Besides, since $\sum_{l=1}^{k}\phi_{l}\phi_{l}^{T}$ follows a Wishart distribution, by Lemma \ref{lemma_c1}, we have $\mathbb{E}[\Vert (\sum_{l=1}^{k}\phi_{l}\phi_{l}^{T})^{-1} \Vert^{4p}]=O(1/k^{4p})$.
Note that $\Vert P_{k} \Vert \leq  \Vert (\sum_{l=1}^{k}\phi_{l}\phi_{l}^{T})^{-1} \Vert/\uline{\beta}$. Then, it holds that $\mathbb{E}[\Vert P_{k} \Vert^{4p}] = O(1/k^{4p}).$

By H\"{o}lder's inequality \cite{ash2014real}, $k^{p}\mathbb{E}[ \Vert P_{k} \Vert^{2p} \Vert \sum_{l=1}^{k}\beta_{l}x_{l} \Vert^{2p}]  \leq k^{p}(\mathbb{E}[ \Vert P_{k} \Vert^{4p} ]\mathbb{E}[ \Vert \sum_{l=1}^{k}\beta_{l}x_{l} \Vert^{4p} ])^{1/2} = O(1)$, which implies that $\mathbb{E}[\Vert \sqrt{k} \tilde{\gamma}_{k}\Vert^{2p}] = O(1).$

By Theorem \ref{mai_the_b} and $\rho(\hat{\delta}_{k}) - \rho(\delta_{y}) = \dot{\rho}(\eta_{k})\tilde{\delta}_{k}$, the boundness of $\vert \dot{\rho}(\eta_{k}) \vert $ implies that $\mathbb{E} [ \vert \sqrt{k} (\rho(\hat{\delta}_{k}) - \rho(\delta_{y}))\vert^{2p}] = O(1).$ 

Then, by (\ref{jioj100}) and the $C_{r}$-inequality \cite{P}, we obtain \begin{align*}
& \mathbb{E}\left[\left\Vert \sqrt{k} \tilde{\theta}_{k} \right\Vert^{2p} \right] \\ \leq & 2^{2p-1} \mathbb{E} \left[  \left\Vert \frac{\sqrt{k}\tilde{\gamma}_{k}}{\rho(\hat{\delta}_{k})} \right\Vert^{2p} + \left\Vert \frac{\sqrt{k}(\rho(\hat{\delta}_{k}) - \rho(\delta_{y}))\gamma}{\rho(\hat{\delta}_{k})\rho(\delta_{y})} \right\Vert^{2p}  \right] = O(1).
\end{align*}
Furthermore, by H\"{o}lder's inequality \cite{ash2014real}, it holds that $\mathbb{E}[\Vert \sqrt{k} \tilde{\theta}_{k} \Vert^{p} ] \leq (\mathbb{E}[\Vert \sqrt{k} \tilde{\theta}_{k} \Vert^{2p} ])^{1/2} = O(1)$,
which indicates (\ref{the_1_21}) holds, i.e., the $L^{p}$ convergence rate is $O(1/k^{p/2})$.

The proof is completed. 

 \subsection{Proof of Proposition \ref{pri_the_uu}}\label{app_2uu}

Since \( A(q) \) and \( B(q) \) are coprime, there exists a specific row in the matrix \( [b_0, b_1, \dots, b_{n_b}] \) that is coprime to \( [1, a_1, \dots, a_{n_a}] \). Without loss of generality, we denote this row vector as \( [b^{*}_{0}, b^{*}_{1}, \ldots, b^{*}_{n_b} ] \in \mathbb{R}^{n_b + 1} \).  
Define the polynomial $B^{*}({ q}) \triangleq b^{*}_{0} + b^{*}_{1}q^{-1} + b^{*}_{ n_b}q^{- n_b}$.
Correspondingly, let \( H^{*}({ q}) \) be the transfer function $H^{*}(q) \triangleq B^{*}(q)/A(q) = \sum_{i=0}^{\infty}h_{i}^{*}q^{-i}$. 
Define the matrix \( \Gamma^{*} \) as
\[
\Gamma^{*} \triangleq
\begin{bmatrix}
h^{*}_{n_b} & h^{*}_{n_b-1} & \cdots & h^{*}_{ n_b+1- n_a} \\
h^{*}_{n_b+1} & h^{*}_{n_b} & \cdots & h^{*}_{n_b+2- n_a} \\
\vdots & \vdots & \ddots & \vdots \\
h^{*}_{n_b + n_a-1} & h^{*}_{ n_b+ n_a-2} & \cdots & h^{*}_{n_b}
\end{bmatrix} \in \mathbb{R}^{ n_a \times  n_b}.
\]  

By a similar argument to Lemma 1 in \cite{zhao2012markov}, it follows that \( \text{rank}(\Gamma^{*}) = n_a \). Additionally, since \( H^{*}(q) \) is constructed as a component of the original transfer function, each row of \( \Gamma^{*} \) corresponds to a row in \( \Gamma \). Consequently, \( \text{rank}(\Gamma) =  n_a \).  

The proof is completed. 

\subsection{Proof of Theorem \ref{thea2}}\label{app_2}

In this proof, we will analyze the properties of the original non-independent and non-stationary sequence $\{s_{k}\}_{k=1}^{\infty}$ by constructing a geometrically strongly mixing and strictly stationary sequence $\{s_{k}^{*}\}_{k=-\infty}^{\infty}$. By examining the difference between $s_{k}^{*}$ and $s_{k}$ and leveraging the favorable characteristics of $\{s_{k}^{*}\}_{k=-\infty}^{\infty}$, the almost sure convergence rate of the proposed algorithm is established by using Lemma \ref{lemma_f3}.  

We define the auxiliary process  
\[
y_{k}^{*} = (\phi_k^{*})^T H(q) + d_{k}^{*}, \quad s_k^{*} = Q(y_k^{*}), \quad k \in \mathbb{Z},
\]  
where $\{\phi_{k}^{*}\}_{k = -\infty}^{\infty}$ and $\{d_{k}^{*}\}_{k = -\infty}^{\infty}$ are two mutually independent and iid sequences. 
For $k \geq 1$, $\phi_{k}^{*} \triangleq \phi_{k}$ and $d_{k}^{*} \triangleq d_{k}$, while for $k \leq 0$, $\phi_k^{*} \sim \mathcal{N}(\mathbf{0}_{n}, H)$ and $d_k^{*} \sim \mathcal{N}(0, \delta_{d}^{2})$. The quantizer $Q(\cdot)$ is defined by (\ref{model_b}).

Now we examine the deviation between the original quantized observation $s_{k}$ and its auxiliary counterpart $s_{k}^{*}$.
Given that $A(q)$ has no poles on or outside the unit circle, there exists a constant $\nu > 0$ such that $\|h_i\| = O(\exp(-\nu i))$.
Define $\delta_k \triangleq (\mathbb{E}[y_k^2])^{1/2}$, which converges to a finite limit $\delta_\infty = (\mathbb{E}[(y_k^*)^2])^{1/2}$.
Let $\tilde{y}_k^* \triangleq y_k^* - y_k = \sum_{i=k}^{\infty} (\phi_{k-i}^*)^T h_i$ and $\tilde{\delta}_k^* \triangleq \mathbb{E}([(\tilde{y}_k^*)^2])^{1/2}$. It follows that $(\tilde{\delta}_{k}^{*})^{2} = \sum_{i=k}^{\infty} h_{i}^{T}Hh_{i}
= O(\exp(-2\nu k))$, which implies that $\tilde{\delta}_{k}^{*} = O(\exp(-\nu k))$.

Let $\varepsilon_k \triangleq \tilde{\delta}_k^* (16\log k)^{1/2}$. For sufficiently large $k$, we have $\varepsilon_k < \delta_k / k^2$, which implies that
\begin{align*}
\mathbb{P}\left( \left\vert y_{k} \right\vert < \varepsilon_{k}\right) & \leq \mathbb{P}\left( \left\vert \frac{y_{k}}{\delta_{k}} \right\vert < \frac{1}{k^2}\right) = \int_{-1/k^2}^{1/k^2} f(x)  dx  \leq \frac{2}{k^2}.
\end{align*}
Besides, since $\tilde{y}_k^*$ is a Gaussian variable, by Lemma \ref{lemma_d1}, it holds that $\mathbb{P}(\vert \tilde{y}_{k}^{*} \vert > \varepsilon_{k}) = \mathbb{P}(\vert \tilde{y}_{k}^{*}/\tilde{\delta}_{k}^{*} \vert > (16\log k)^{1/2}) = O(1/k^{2})$. 
Since a discrepancy $s_k \neq s_k^*$ arises if either $y_k$ lies within a distance of $\varepsilon_k$ from a threshold, or the deviation $\tilde{y}_k^* = y_k^* - y_k$ exceeds $\varepsilon_k$, one can get $\mathbb{P}(s_{k} \neq s_{k}^{*}) \leq \mathbb{P}( \bigcup_{j=1}^{m} \{ \vert y_{k} - C_{j}\vert < \varepsilon_{k}\}) + \mathbb{P}( \vert \tilde{y}_{k}^{*} \vert > \varepsilon_{k}) \leq m\mathbb{P}( \vert y_{k} \vert < \varepsilon_{k} ) + \mathbb{P}(\vert \tilde{y}_{k}^{*} \vert > \varepsilon_{k} )$. It follows that
\begin{align}\label{YTREEs}
\mathbb{P}\left( s_{k} \neq s_{k}^{*} \right) = O\left(\frac{1}{k^2}\right).
\end{align}

Now we examine the strictly stationary and geometrically strongly mixing properties of the auxiliary sequence $\{s_{k}^{*}\}_{k=-\infty}^{\infty}$.
Since \( A(q) \) has no poles on or outside the unit circle and $\{\phi_{k}^{*}\}_{k=-\infty}^{\infty}$ is iid Gaussian, $\{y_{k}^{*}\}_{k=-\infty}^{\infty}$ is strictly stationary.
Besides, since $s_k^*$ depends only on the quantization interval into which $y_k^*$ falls, $\{s_k^*\}$ is also strictly stationary.
We now analyze the geometrically strongly mixing property.
First, we decorrelate $\phi_k^*$ by defining $u_{k}^{*} = [u_{k,1}^*,u_{k,2}^*,\ldots,u_{k,n}^*]^{T}  \triangleq H^{-1/2} \phi_k^* \sim \mathcal{N}(\mathbf{0}_n, I_n)$. 
Then, following Lemma 1 and Theorem 1' in \cite{MOKKADEM1988309}, for the system $A(q)z_{k,i}^{*} = (H^{1/2}B(q))_{i}u_{k,i}^{*}$, where $(H^{1/2}B(q))_{i}$ denotes its $i$-th component, there exists a Markov chain $\{m_{k,i}\in \mathbb{R}^{\max\{n_a,n_b+1\}}\}_{k=-\infty}^{\infty} $ such that
\[
m_{k,i} = D_{1,i}m_{k-1,i} + D_{2,i}u_{k,i}^{*}, \quad z_{k,i}^{*} = D_{3,i}m_{k,i},
\]
for suitable matrices or vectors $D_{1,i}, D_{2,i}$ and $D_{3,i}$.
Moreover, $\{m_{k,i} \}_{k=-\infty}^{\infty}$ is geometrically strongly mixing, i.e., there exists a positive constant $\varkappa_{i}< 1$ such that its $\alpha$-mixing coefficient defined by (\ref{alphamix}) satisfies $\alpha(\tau) = O(\varkappa_i^\tau)$.
Additionally, due to the independence of $u_{k,i}^*$ and $u_{k,j}^*$ for $i \neq j$, one can get $m_{k,i}^*$ and $m_{k,j}^*$ are also independent.
Thus, there exists a positive constant $\varkappa < 1$ such that
\[
\sup_{t \in \mathbb{Z}} \sup_{A \in \mathcal{G}_{-\infty}^{t}, B \in \mathcal{G}_{t + \tau}^{\infty}} \left\vert \mathbb{P}(A\cap B) - \mathbb{P}(A)\mathbb{P}(B) \right\vert
=  O(\varkappa^\tau),
\]
where $ \mathcal{G}_{-\infty}^{t} \triangleq \sigma\{(m_{k,1}, m_{k,2} \ldots, m_{k,n}): k \leq t\}$ and $\mathcal{G}_{t+\tau}^{\infty} \triangleq \sigma\{(m_{k,1}, m_{k,2} \ldots, m_{k,n}): k \geq t+\tau\}$.
Besides, for $1 \leq i \leq n$, since $D_{2,i} \neq \mathbf{0}_{\max\{n_a, n_b+1\}}$, $u_{k,i}^{*}$ is measurable with respect to $ m_{k,i} $ and $m_{k-1,i}$.
Since $\phi_k^{*}$ is measurable with respect to $u_k^{*}$, $z_{k,i}^{*}$ is measurable with respect to $m_{k,i}$, $y_k^{*} = \sum_{i=1}^{n} z_{k,i}^{*} + d_k^{*}$, and $d_k^{*}$ is independent of $m_{k,i}$ for $1 \leq i \leq n$, one can get 
\begin{align}\label{yy6ees2}
\sup_{t \in \mathbb{Z}} \sup_{ 
\substack{A \in \mathcal{F}_{-\infty}^t, \\ B \in \mathcal{F}_{t+\tau+1}^\infty} }\left\vert \mathbb{P}(A\cap B) - \mathbb{P}(A)\mathbb{P}(B) \right\vert
=  O(\varkappa^\tau),
\end{align}
where $ \mathcal{F}_{-\infty}^t \triangleq \sigma\{(\phi_{k}^{*},d_{k}^{*}): k \leq t\}$ and $\mathcal{F}_{t+\tau+1}^{\infty} \triangleq \sigma \{(\phi_{k}^{*}, y_{k}^{*}): k \geq t+\tau+1 \}$.

Based on (\ref{YTREEs}) and (\ref{yy6ees2}), the following proof will be divided into the following two parts.

$\mathbf{Part}$ $\mathbf{1:}$ Establish the almost sure convergence rate of estimating $\delta_{\infty}$.

Given the strict stationarity of $\{s_k^{*}\}_{k=-\infty}^{\infty}$, $\sigma\{s_{k}^{*}: k \leq t\} \subseteq \mathcal{F}_{-\infty}^{t}$, and $\sigma\{s_{k}^{*}: k \geq t+ \tau+1\} \subseteq  \mathcal{F}_{t+\tau+1}^\infty$, for $j=1,2,\ldots,m$, the auxiliary sequence $\{w_{k,j}^{*} \triangleq I_{\{s_{k}^{*} < j \}} - F(C_{j}/\delta_{\infty})\}_{k=-\infty}^{\infty}$ is also strictly stationary and geometrically strongly mixing. 
Therefore, its càdlàg inverse function of the $\alpha$-mixing coefficient satisfies $\alpha_{w,j}^{-1}(u) = O(-\log u ).$
By $\mathbb{Q}_{w,j}(u) \triangleq \inf \{ \tau \geq 0 : \mathbb{P}(\vert w_{0,j}^{*} \vert > \tau) \leq u \} < \vert w_{0,j}^{*} \vert < 2$, 
\[
\int_0^1 \alpha_{w,j}^{-1}(u) \mathbb{Q}_{w,j}^2(u) du = O\left(\int_0^1 -\log u du\right) < \infty.
\]
Since $\{w_{k,j}^{*}\}_{k=-\infty}^{\infty}$ has zero mean, by Lemma \ref{lemma_f3}, we have
\[
\left\vert\sum_{l=1}^{k} w_{l,j}^{*} \right\vert =  O\left(\sqrt{k\log\log k}\right), \quad \text{a.s.} 
\]
Besides, since $\mathbb{P}( \vert s_{k} - s_{k}^{*} \vert \geq 1/k^{2} ) = \mathbb{P}( s_{k} \neq s_{k}^{*})$,  by (\ref{YTREEs}) and the Borel-Cantelli Lemma \cite{ash2014real}, we have
\[
\sum_{l=1}^{k} \left\vert  I_{\{s_{l} < j \}} - I_{\{s_{l}^{*}<j\}} \right\vert = O\left( \sum_{l=1}^{k} \frac{1}{l^2} \right) = O(1), \quad \text{a.s.}
\]
From $\sum_{i=0}^{j-1}S_{k}^{i} - F(C_{j}/\delta_{\infty}) = (\sum_{l=1}^{k} w_{l,j}^{*} +  I_{\{s_{l} < j \}} - I_{\{s_{l}^{*}<j\}})/k$, one can get
\[
 \left\vert\sum_{i=0}^{j-1} S_{k}^{i} - F\left(\frac{C_{j}}{\delta_{\infty}}\right)\right\vert 
= O\left(\sqrt{\frac{\log\log k}{k}}\right), \quad  \text{a.s.} 
\]

Without loss of generality, we suppose that $C_{j} \neq 0$.
Then, by (\ref{ytrt2}) and the mean value theorem \cite{zorich2016mathematical}, 
for $\hat{\delta}_{k}^{j}$ given by the ML-type algorithm, where $j = 1,2,\ldots,m$, there exists a point $\psi_{k,j}$ between $\sum_{i=0}^{j-1}S_{k}^{i}$ and $F(C_{j}/\delta_{\infty})$ such that $\hat{\delta}_{k}^{j} - \delta_{\infty} = \delta_{\infty}g(\psi_{k,j})(F(C_{j}/\delta_{\infty})-\sum_{i=0}^{j-1} S_{k}^{i})/F^{-1}(\sum_{i=0}^{j-1} S_{k}^{i})$.
Similarly to Section V-C, one can get
\[
\left\vert \hat{\delta}_{k}^{j} - \delta_{\infty} \right\vert= O\left(\sqrt{\frac{\log\log k}{k}}\right), \quad  \text{a.s.},
\]
where the above formula can be proven similarly when $C_{j} = 0$.
Since the positive constant $c$ is sufficiently small, we assume $\delta_{\infty} \in (c,1/c)$.
Then, by (\ref{www23456}), (\ref{boundmu}) and Lemma \ref{lemma_b1}, for $\hat{\delta}_{k}$ given by the ML-type algorithm, we have
\begin{align}\label{yyy2765491}
\left\vert\hat{\delta}_{k} - \delta_{\infty} \right\vert
& = O\left(\sqrt{\frac{\log\log k}{k}}\right), \quad  \text{a.s.}
\end{align}

$\mathbf{Part}$ $\mathbf{2:}$ Establish the almost sure convergence rates of the DM-type algorithm.

Define $\varphi_k^{*} \triangleq [(\phi_k^*)^T, (\phi_{k-1}^*)^T, \ldots, (\phi_{k-\kappa}^*)^T]^T \in \mathbb{R}^{(\kappa + 1)n}$ for all $k \in \mathbb{R}$. 
Similarly to Proposition \ref{pri_the_b}, one can get
\begin{align}\label{e8e5}
 \mathbb{E}\left[s_{k}^{*}\varphi_{k}^{*}\right] =  \rho(\delta_{\infty})\mathbb{E}\left[\varphi_{k}^{*}(\varphi_{k}^{*})^{T}\right]h  = \rho(\delta_{\infty}) \mathbf{H} h,
\end{align}
where $\mathbf{H} = \text{diag}(H,H,\ldots,H) \in \mathbb{R}^{n(\kappa+1) \times n(\kappa+1)}$.

We introduce an another auxiliary variable defined as $x_{k}^{*} \triangleq s_{k}^{*}\varphi_{k}^{*} -  \rho(\delta_{\infty}) \varphi_{k}^{*}(\varphi_{k}^{*})^{T}h$.
Since $\sigma\{x_{k}^{*}: k \leq t\} \subseteq \mathcal{F}_{-\infty}^{t}$, $\sigma\{x_{k}^{*}: k \geq t+\tau+\kappa+1\} \subseteq \mathcal{F}_{t+\tau+\kappa+1}^\infty$, $\{s_k^{*}\}_{k=-\infty}^{\infty}$ and $\{\varphi_{k}^{*}\}_{k=-\infty}^{\infty}$ are both strictly stationary, we know $\{x_{k}^{*}\}_{k=-\infty}^{\infty}$ is also both strictly stationary and geometrically strongly mixing.
Thus, its càdlàg inverse function of the $\alpha$-mixing coefficient also satisfies $\alpha_{x}^{-1}(u) = O(-\log u )$.
Note that $\mathbb{P}( \Vert x_{0}^{*}\Vert  > \tau) \leq \mathbb{P}(\vert s_0^{*} \vert > (\tau/2)^{1/2}) + \mathbb{P}(\vert \rho(\delta_{\infty})\varphi_{0}^{T}h \vert > (\tau/2)^{1/2}) + 2 \mathbb{P}( \Vert \varphi_{0} \Vert > (\tau/2)^{1/2})$.
Since $\varphi_{0}^{*}$ is a Gaussian variable and $s_{0}^{*}$ is uniformly bounded, by Lemma \ref{lemma_d1}, it holds that $\mathbb{Q}_{x}(u) \triangleq \inf \{ \tau \geq 0 : \mathbb{P}(\Vert x_{0}^{*}\Vert > \tau) \leq u \} = O(- \log u)$. 
By (\ref{e8e5}) and Lemma \ref{lemma_f3}, one can get
\[
\left\Vert\sum_{l=1}^{k} x_{l}^{*} \right\Vert = O\left(\sqrt{k\log\log k}\right), \quad \text{a.s.} 
\]
Besides, for $k \geq \kappa + 1$, since $\mathbb{P}( \vert s_{k}\varphi_{k} - s_{k}^{*}\varphi_{k}^{*} \vert > 1/k^{2} ) = \mathbb{P}( \vert (s_{k} - s_{k}^{*})\varphi_{k} \vert > 1/k^{2} ) \leq \mathbb{P}( s_{k} \neq s_{k}^{*})$, by (\ref{YTREEs}) and the Borel-Cantelli Lemma \cite{ash2014real}, 
\[
\sum_{l=1}^{k} \left\Vert  s_{l}\varphi_{l}- s_{l}^{*}\varphi_{l}^{*} \right\Vert = O(1), \quad \text{a.s.}
\]
Hence,
\[
\frac{1}{k}\left\Vert \sum_{l=1}^{k}  s_{l}\varphi_{l}- \rho(\delta_{\infty})\varphi_{l}\varphi_{l}^{T}h \right\Vert  =  O\left(\sqrt{\frac{\log\log k}{k}}\right), \quad  \text{a.s.} 
\]

Note that $\sum_{l=1}^{k}\varphi_{l}\varphi_{l}^{T}$ can be rewriten as $\sum_{l=1}^{k}\varphi_{l}\varphi_{l}^{T} = \sum_{u=1}^{\kappa+1} \sum_{v=0}^{\lceil k/(\kappa+1)\rceil} \varphi_{u+v(\kappa+1)}\varphi_{u+v(\kappa+1)}^{T}$, where we stipulate that $\varphi_{w} = 0$ for $w > k$. 
Then, since for $u = 1,2,\ldots,\kappa+1$, each sequence $\{ \varphi_{u+v(\kappa+1)}\varphi_{u+v(\kappa+1)}^{T} \}_{v = 1}^{\lceil k/(\kappa+1)\rceil}$ is an iid sequence with finite variance, the law of the iterated logarithm \cite{chung2000course} implies that $\sum_{l=1}^{k}\varphi_{l}\varphi_{l}^{T}/k$ converges to $\mathbf{H}$ almost surely, which indicates that $\Vert k (\sum_{l=1}^{k}\varphi_{l}\varphi_{l}^{T})^{-1}  \Vert$ converges to $\Vert \mathbf{H}^{-1}\Vert$ almost surely.
Thus, for $P_{k}$ given by the WLS-type algorithm, by $\Vert P_{k} \Vert \leq \Vert (\sum_{l=1}^{k}\varphi_{l}\varphi_{l}^{T})^{-1} \Vert$, we have (\ref{ppppi}) holds.
Then, for $\hat{\gamma}_{k}$ given by the WLS-type algorithm, similarly to Section V-C, 
\[
\left\vert \hat{\gamma}_{k} - \rho(\delta_{\infty})\theta \right\vert =O\left(\sqrt{\frac{\log\log k}{k}}\right), \quad  \text{a.s.}
\]
Together with (\ref{yyy2765491}), similarly to the proof of Theorem \ref{mai_the_e}, 
\[
\left\Vert \hat{h}_{k} - h \right\Vert =O\left(\sqrt{\frac{\log\log k}{k}}\right), \quad \text{a.s.}
\]
Finally, by Proposition \ref{pri_the_uu}, similarly with the proof of Theorem 2 in \cite{song2018recursive}, we have (\ref{aa_pp_2}) holds.

The proof is completed.

\section{Numerical examples}\label{sec_d}
This section will illustrate the theoretical results with simulation examples.

$\mathbf{Example}$ $\mathbf{1}$: Consider a linear system with the quantized observation as:
\[
y_{k} = \phi_{k}^{T}\theta + d_{k}, \quad s_{k} = I_{\{ y_{k} \leq 0 \} } + I_{\{ y_{k} \leq 1 \} },
\]
where the unknown parameter is $\theta = [0.2,-0.1,0.5]^{T}$; the system noise $d_{k}$ follows $\mathcal{N}(0,0.25)$; the system input $\phi_{k}$ follows $\mathcal{N}(\mathbf{0}_{3},H)$, where
$$
H =
\begin{bmatrix}
2.5 &-0.6 &-0.4  \\
-0.6 &2 &0.6 \\
-0.4 &0.6 &1.5 \\
\end{bmatrix}.
$$

\begin{figure}[h]
\centering
\noindent\includegraphics[width=0.48\textwidth]{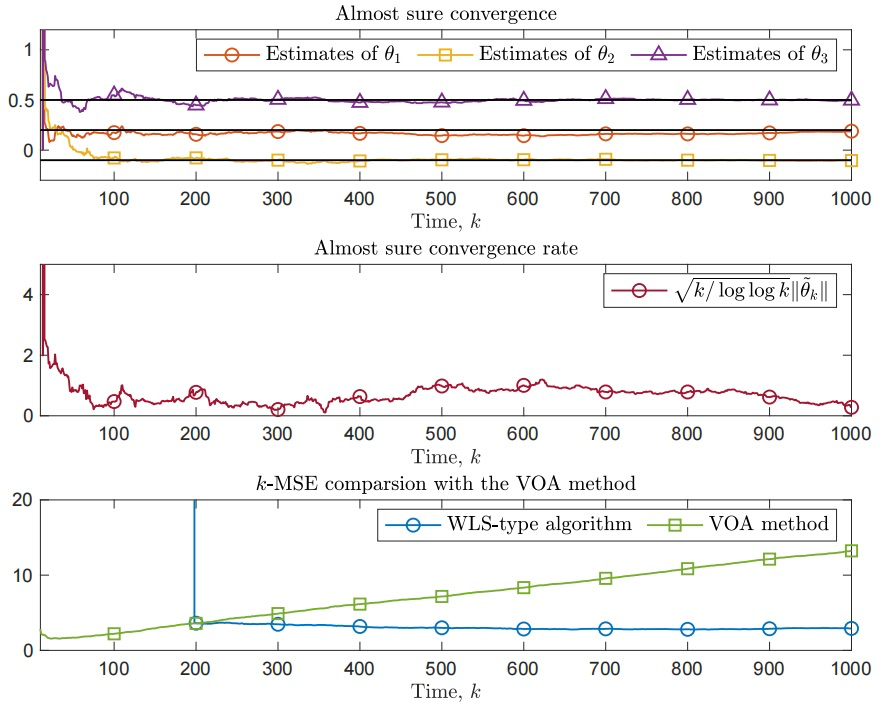}
\caption{Convergence of the WLS-type algorithm for estimating $\theta$}
\label{fig_b}
\end{figure}

Here we apply the ML-type algorithm and the WLS-type algorithm to give the estimates, initializing with $P_{0} = I_{3}/10$ and $\hat{\gamma} = [0,0,0]^{T}$, and using weight coefficients $\beta_{k} \equiv 1$ and projection operator parameters $c = 10^{-6}$ and $c^{*} = 10^{-6}$. 
Figure \ref{fig_b} illustrates the estimates given by the WLS-type algorithm, demonstrating the boundedness of the trajectory of $\sqrt{k\log\log k} \Vert \tilde{\theta}_{k} \Vert$, which indicates almost sure convergence with a convergence rate of $O(\sqrt{\log\log k/k})$.

Besides, we conducted 250 Monte Carlo simulations to compare the mean square error (MSE) of the WLS-type algorithm with the variationally optimal approximation (VOA) method introduced in \cite{risuleo2019identification}.
As illustrated in Figure \ref{fig_b}, the trajectory of the $\mathbb{E}[\tilde{\theta}_{k}^{T}\tilde{\theta}_{k}]$ exhibits a convergence rate of $O(1/k)$.
Furthermore, Fig. \ref{fig_b} demonstrates that while the MSE of the WLS-type algorithm is initially higher than that of the VOA method, it surpasses the VOA in performance as iterations proceed.
This result is consistent with the analysis in \cite{risuleo2019identification}, where it was shown that the VOA method might fail to achieve an MSE convergence rate of $O(1/k)$.

\begin{figure}[h]
\centering
\noindent\includegraphics[width=0.48\textwidth]{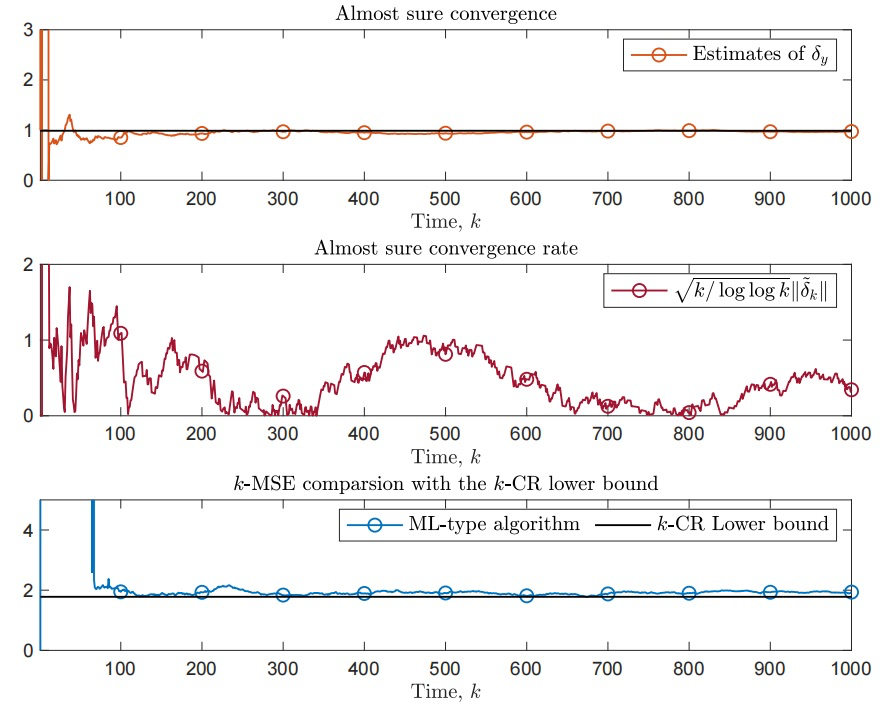}
\caption{Convergence of the ML-type algorithm for estimating $\delta_y$}
\label{fig_c}
\end{figure}

In addition, the convergence of the ML-type algorithm is shown in Fig. \ref{fig_c}, demonstrating $\hat{\delta}_{k}$ converges to $\delta_{y}$ almost surely with the convergence rate of $O(\sqrt{\log \log k /k})$ and is an asymptotically efficient
estimate for estimating the variance of Gaussian variables based on quantized observations.

$\mathbf{Example}$ $\mathbf{2}$: Consider an OE system with the quantized observation as:
\begin{align*}
& y_{k} = \phi_{k}(b_{0} +  b_{1}z + b_{2}z^{2})/(1+a_{1}z) + d_{k}, \\
& s_{k} = I_{\{y_{k} \leq 1.5 \} },
\end{align*}
where the unknown parameter is $\theta^{*} = [a_{1},b_{0},b_{1},b_{2}]^{T} = [0.2, 1,-0.2,0.6]^{T}$; the system noise $d_{k}$ follows $\mathcal{N}(0,1)$; the system input $\phi_{k}$ follows $\mathcal{N}(0,1)$.

\begin{figure}[h]
\centering
\noindent\includegraphics[width=0.48\textwidth]{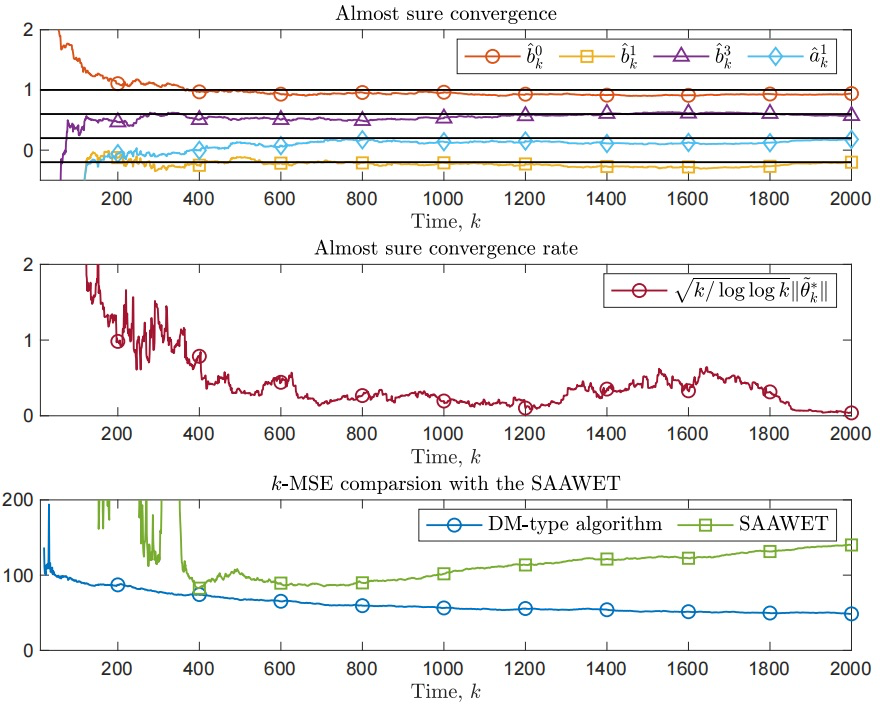}
\caption{Convergence of the DM-type algorithm for estimating $\theta^{*}$}
\label{fig_d}
\end{figure}

We employ DM-type algorithms to obtain parameter estimates. 
Fig. \ref{fig_d} illustrates the convergence of the DM-type algorithm, showing almost sure convergence with a rate of $O(\sqrt{\log \log k/k})$. 
To evaluate its estimation accuracy, we perform 250 simulations comparing it with the stochastic approximation algorithm with expanding truncations (SAAWET) from \cite{song2018recursive}. 
The results in Fig. \ref{fig_d} indicate that the DM-type algorithm, leveraging second-order information, achieves better MSE near the convergence point. 
In contrast, while SAAWET can converge to a minimum using single-sample updates, its reliance on first-order information prevents it from attaining an MSE convergence rate of $O(1/k)$.

\section{Conclusion} \label{sec_g}
This paper introduces a WLS-type criterion for quantized identification problems and proposes a corresponding unified recursive identification algorithm. 
This identification algorithm can adapt to noisy and noise-free scenarios under fixed thresholds and Gaussian inputs.
Under mild conditions, the convergence properties of this identification algorithm are established.
Moreover, this algorithm offers an asymptotically efficient estimation of the variance of Gaussian variables based on quantized observations.
Furthermore, the proposed methods and results are extended to OE systems.

Here we give some topics for future research. 
Firstly, the design of the weight coefficients is left as an open question. 
Second, how can we assess the explicit expression of the asymptotic variance?
Thirdly, how can we remove the Gaussian restriction of the input and noise?

\vspace{-4 em}
\begin{IEEEbiography}[{\includegraphics[width=1in, height=1.25in, clip, keepaspectratio]{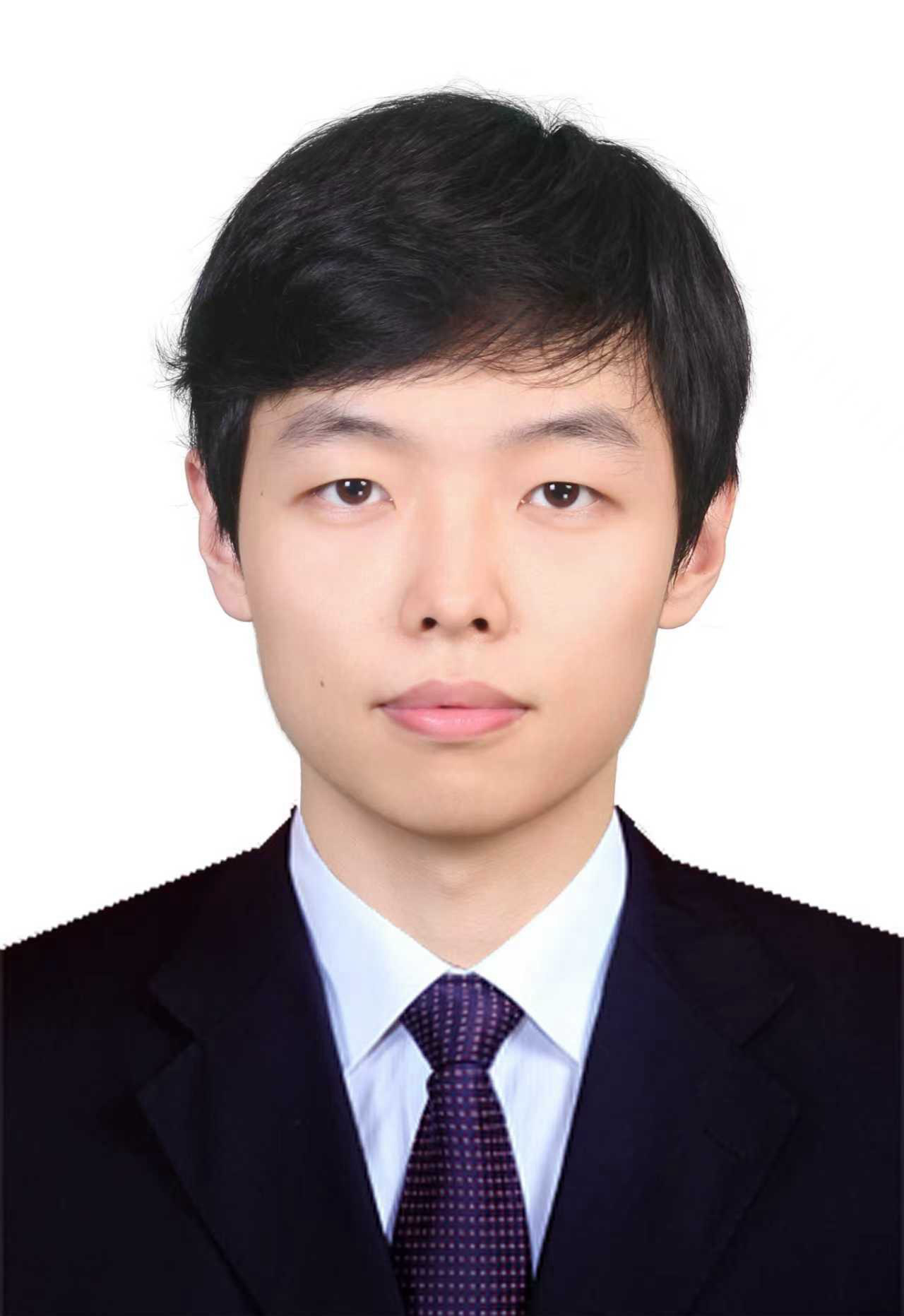}}]{Xingrui Liu} (Student Member, IEEE) received the B.S. degree in statistics from the School of Mathematics, Shandong University, Jinan, China, in 2022. He is working toward a Ph.D. degree in system science at the Academy of Mathematics and Systems Science, Chinese Academy of Science, Beijing, China. 

His research interests include identification, control, and communication in quantized systems.
\end{IEEEbiography}
\vspace{-4 em}
\begin{IEEEbiography}[{\includegraphics[width=1in, height=1.25in, clip, keepaspectratio]{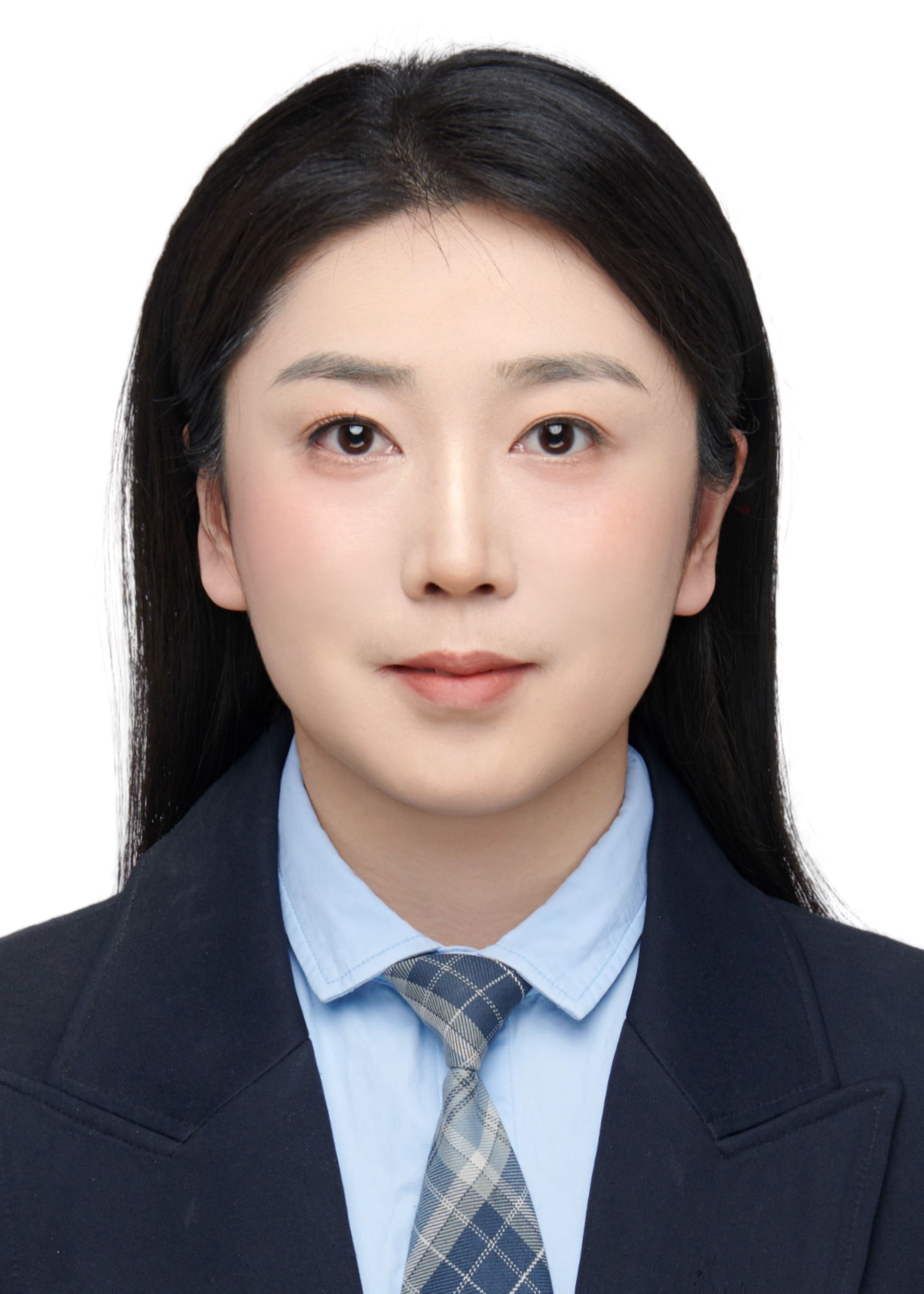}}]{Ying Wang} (Member, IEEE) received the B.S. degree in mathematics from Wuhan University, Wuhan, China, in 2017, and the Ph.D. degree in system theory from the Academy of Mathematics and Systems Science (AMSS), Chinese Academy of Sciences (CAS), Beijing, China, in 2022. 

She is currently a postdoctoral fellow at the AMSS, CAS, and Division of Decision and Control Systems, KTH Royal Institute of Technology, Stockholm, Sweden. Her research interests include parameter identification and adaptive control of quantized systems. She was a recipient of the Special Prize of the Presidential Scholarship of Chinese Academy of Sciences in 2022, and the Outstanding Doctoral Dissertation Award of the Chinese Association of Automation in 2024.

\end{IEEEbiography}
\vspace{-4 em}
\begin{IEEEbiography}[{\includegraphics[width=1in, height=1.25in, clip, keepaspectratio]{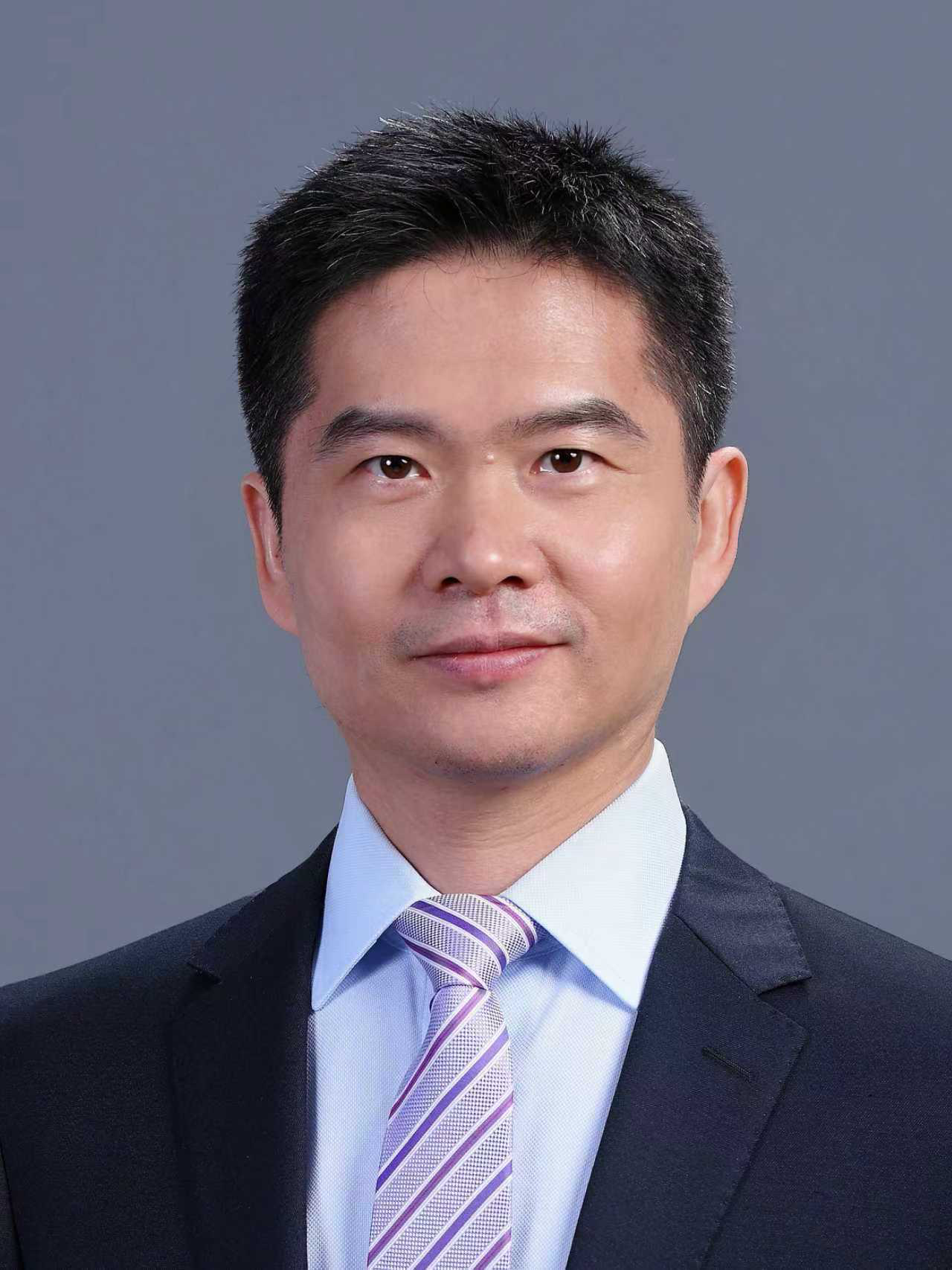}}]{Yanlong Zhao} (Senior Member, IEEE) received the B.S. degree in mathematics from Shandong University, Jinan, China, in 2002, and the Ph.D. degree in systems theory from the Academy of Mathematics and Systems Science (AMSS), Chinese Academy of Sciences (CAS), Beijing, China, in 2007. 

Since 2007, he has been with the AMSS, CAS, where he is currently a full Professor and a Vice Director of State Key Laboratory of Mathematical Sciences. 
His research interests include the identification and control of quantized systems, information theory and modelling of financial systems. 

Yanlong Zhao was awarded the Second Prize of the State Natural Science Award of China in 2015. He has been a Deputy Editor-in-Chief of {\em Journal of Systems Science and Complexity}, an Associate Editor for {\em Automatica}, {\em SIAM Journal on Control and Optimization}, and {\em IEEE Transactions on Systems, Man and Cybernetics: Systems}. He served as a Vice-President of Asian Control Association and a Vice-President of IEEE CSS Beijing Chapter, and is now a Vice-President and Fellow of the Chinese Association of Automation (CAA) and the Chair of Technical Committee on Control Theory, CAA.
\end{IEEEbiography}
\vspace{-0.5 em}
\end{document}